\documentclass[letterpaper, 10pt, conference]{ieeeconf}
\IEEEoverridecommandlockouts 
\usepackage{amsmath,amssymb,url}
\usepackage{graphicx,subfigure}

\newcommand{\norm}[1]{\ensuremath{\left\| #1 \right\|}}

\newcommand{\bracket}[1]{\ensuremath{\left[ #1 \right]}}

\newcommand{\parenth}[1]{\ensuremath{\left( #1 \right)}}

\newcommand{\refeqn}[1]{(\ref{eqn:#1})}
\newcommand{\reffig}[1]{Fig. \ref{fig:#1}}
\newcommand{\tr}[1]{\mbox{tr}\ensuremath{\negthickspace\bracket{#1}}}

\newcommand{\SO}{\ensuremath{\mathsf{SO(3)}}}
\newcommand{\T}{\ensuremath{\mathsf{T}}}

\newcommand{\so}{\ensuremath{\mathfrak{so}(3)}}

\renewcommand{\Re}{\ensuremath{\mathbb{R}}}
\newcommand{\Sph}{\ensuremath{\mathsf{S}}}

\newcommand{\D}{\ensuremath{\mathbf{D}}}

\title{\LARGE \bf
Stable Manifolds of Saddle Points\\ for Pendulum Dynamics on $\Sph^2 $ and $\SO$}

\author{Taeyoung Lee\authorrefmark{1}, Melvin Leok\authorrefmark{2}, and N. Harris McClamroch%
\thanks{Taeyoung Lee, Mechanical and Aerospace Engineering, Florida Institute of Technology, Melbourne, FL 39201 {\tt taeyoung@fit.edu}}%
\thanks{Melvin Leok, Mathematics, University of California at San Diego, La Jolla, CA 92093 {\tt mleok@math.ucsd.edu}}%
\thanks{N. Harris McClamroch, Aerospace Engineering, University of Michigan, Ann Arbor, MI 48109 {\tt
nhm@umich.edu}}%
\thanks{\textsuperscript{\footnotesize\ensuremath{*}}This research has been supported in part by NSF under grants CMMI-1029551.}
\thanks{\textsuperscript{\footnotesize\ensuremath{\dagger}}This research has been supported in part by NSF under grants DMS-0726263, DMS-1001521, DMS-1010687, and CMMI-1029445.}
}

\begin{document}
\allowdisplaybreaks
\maketitle \thispagestyle{empty} \pagestyle{empty}

\begin{abstract}
Attitude control systems naturally evolve on nonlinear configuration spaces, such as $\Sph^2$ and $\SO$. The nontrivial topological properties of these configuration spaces result in interesting and complicated nonlinear dynamics when studying the corresponding closed loop attitude control systems. In this paper, we review some global analysis and simulation techniques that allow us to describe the global nonlinear stable manifolds of the hyperbolic equilibria of these closed loop systems. A deeper understanding of these invariant manifold structures are critical to understanding the global stabilization properties of closed loop control systems on nonlinear spaces, and these global analysis techniques are applicable to a broad range of problems on nonlinear configuration manifolds.
\end{abstract}

\section{Introduction}

Global nonlinear dynamics of various classes of closed loop attitude control systems have been studied in recent years.   An overview of results on attitude control of a rotating rigid body is given in ~\cite{ChSaMcCSM11}.   Closely related results on attitude control of a spherical pendulum (with attitude an element of the two-sphere $\Sph^2$) and of a 3D pendulum (with attitude an element of the special orthogonal group $\SO$) are given in ~\cite{ChMcIJRNC07,ChMcBeAut08,ChaMcCITAC09}.    These and other similar publications address the global closed dynamics of smooth vector fields.   Assuming that the closed loop vector field has an asymptotically stable equilibrium, as desired in attitude stabilization problems, additional hyperbolic equilibria necessarily exist.   The domain of attraction of the asymptotically stable equilibrium is contained in the complement of the union of the stable manifolds of the hyperbolic equilibria. These geometric factors motivate the current paper, in which new analytical and computational results on the stable manifolds of the hyperbolic equilibria are obtained.   

To make the development concrete, the presentation is built around two specific closed loop vector fields: one for the attitude dynamics of a spherical pendulum and one for the attitude dynamics of a 3D pendulum.  In analyzing these two cases, we introduce new analytical and computational tools that are broadly applicable to studying the geometry of more general attitude control systems.  

\section{Spherical Pendulum}

A spherical pendulum is composed of a mass $m$ connected to a frictionless pivot by a massless link of length $l$.
It is acts under uniform gravity, and it is subject to a control moment $u$.  The configuration of a spherical pendulum is described by a unit-vector $q\in\Re^3$, representing the direction of the link with respect to a reference frame.

Therefore, the configuration space is the two-sphere $\Sph^2=\{q\in\Re^3\,|\, q\cdot q =1\}$. The tangent space of the two-sphere at $q$, namely $\T_q\Sph^2$, is the two-dimensional plane tangent to the unit sphere at $q$, and it is identified with $\T_q\Sph^2\simeq\{\omega\in\Re^3\,|\, q\cdot \omega =0\}$, using the following kinematics equation:
\begin{align*}
\dot q= \omega\times q,
\end{align*}
where the vector $\omega\in\Re^3$ represents the angular velocity of the link. The equation of motion is given by
\begin{gather*}
\dot\omega = \frac{g}{l}q\times e_3 + \frac{1}{ml^2} u,
\end{gather*}
where the constant $g$ is the gravitational acceleration, and the vector $e_3=[0,0,1]\in\Re^3$ denotes the unit vector along the direction of gravity. The control moment at the pivot is denoted by $u\in\Re^3$.

\subsection{Control System}

Several proportional-derivative (PD) type control systems have been developed on $\Sph^2$ in a coordinate-free fashion~\cite{BulMurN95,BulLew05}. Here, we summarize a control system that stabilizes a spherical pendulum to a fixed desired direction $q_d\in\Sph^2$. 

Consider an error function on $\Sph^2$, representing the projected distance from the direction $q$ to the desired direction $q_d$, given by
\begin{align*}
\Psi(q,q_d)=1-q\cdot q_d.
\end{align*}
The derivative of $\Psi$ with respect to $q$ along the direction $\delta q =\xi\times q$, where $\xi\in\Re^3$ and $\xi\cdot q=0$, is given by
\begin{align*}
\D_q \Psi(q,q_d)\cdot\delta q = -(\xi\times q)\cdot q_d = (q_d\times q)\cdot \xi.
\end{align*}
For positive constants $k_q,k_\omega$, the control input is chosen as:
\begin{align*}
u = ml^2(- k_\omega \omega - k_q q_d\times q -\frac{g}{l}q\times e_3).
\end{align*}
The corresponding closed loop dynamics are given by
\begin{gather}
\dot \omega = -k_\omega\omega - k_q q_d\times q,\label{eqn:dotw}\\
\dot q= \omega\times q.\label{eqn:dotq}
\end{gather}

This yields two equilibrium solutions: (i) the desired equilibrium $(q,\omega)=(q_d,0)$; (ii) additionally, there exists another equilibrium $(-q_d,0)$ at the antipodal point on the two-sphere. 

It can be shown that the desired equilibrium is asymptotically stable by using the following Lyapunov function:
\begin{align*}
\mathcal{V} = \frac{1}{2} \omega\cdot\omega + k_q \Psi(q,q_d).
\end{align*}
In this paper, we analyze the local stability of each equilibrium by linearizing the closed loop dynamics to study the equilibrium structures more explicitly. In particular, we develop a coordinate-free form of the linearized dynamics of \refeqn{dotw}, \refeqn{dotq}, in the following section.

\subsection{Linearization}

A variation of a curve $q(t)$ on $\Sph^2$ is a family of curves $q^\epsilon(t)$ parameterized by $\epsilon\in\Re$, satisfying several properties~\cite{BulLew05}. It cannot be simply written as $q^\epsilon(t)=q(t)+\epsilon\delta q(t)$ for $\delta q(t)$ in $\Re^3$, since in general, this does not guarantee that $q^\epsilon(t)$ lies in $\Sph^2$. In~\cite{LeeLeoIJNME08}, an expression for a variation on $\Sph^2$ is given in terms of the exponential map as follows:
\begin{align}
q^\epsilon(t) = \exp(\epsilon \hat\xi(t))q(t),\label{eqn:qe}
\end{align}
for a curve $\xi(t)$ in $\Re^3$ satisfying $\xi(t)\cdot q(t)=0$ for all $t$. The \textit{hat map} $\hat\cdot:\Re^3\rightarrow\so$ is defined by the condition that $\hat x y =x\times y$ for any $x,y\in\Re^3$. The resulting infinitesimal variation is given by 
\begin{align}
\delta q (t) = \frac{d}{d\epsilon}\bigg|_{\epsilon=0} q^\epsilon(t) = \xi(t)\times q(t). \label{eqn:delq}
\end{align}
The variation of the angular velocity can be written as
\begin{align}
\omega^\epsilon(t) = \omega(t) + \epsilon \delta \omega(t),\label{eqn:we}
\end{align}
for a curve $\delta w(t)$ in $\Re^3$ satisfying $q(t)\cdot w(t)=0$ for all $t$. Hereafter, we do not write the dependency on time $t$ explicitly. 

The time-derivative of $\delta q$ can be obtained either from \refeqn{delq} or by substituting \refeqn{qe}, \refeqn{we} into \refeqn{dotq}, and considering the first order terms of $\epsilon$. In either case, we have
\begin{align*}
\delta\dot q= \dot\xi\times q + \xi\times(\omega\times q) = \delta\omega\times q + \omega\times(\xi\times q).
\end{align*}
Using the vector cross product identity $a\times (b\times c)=(a\cdot c)b-(a\cdot b) c$ for any $a,b,c\in\Re^3$, this can be written as
\begin{gather*}
\dot\xi\times q + (\xi\cdot q)w -(\xi\cdot\omega) q = \delta\omega\times q +(\omega\cdot q)\xi -(\omega\cdot\xi)q.
\end{gather*}
Since $\xi\cdot q=0$, $\omega\cdot q=0$, this reduces to
\begin{gather*}
\dot\xi\times q = \delta\omega\times q.
\end{gather*}
Since both sides of the above equation are perpendicular to $q$, this is equivalent to $q\times(\dot\xi\times q) = q\times(\delta\omega\times q)$, which yields
\begin{gather*}
\dot \xi - (q\cdot\dot\xi) q = q\times(\delta\omega\times q).
\end{gather*}
Since $\xi\cdot q =0$, we have $\dot\xi\cdot q +\xi\cdot\dot q=0$. Using this, the above equation can be rewritten as
\begin{align}
\dot \xi & = -(\xi\cdot(\omega\times q))q+ q\times(\delta\omega\times q)\nonumber\\
& = (qq^T\hat\omega) \xi +(I-qq^T)\delta\omega.\label{eqn:dotxi}
\end{align}
This corresponds to the linearized equation of motion for \refeqn{dotq}. Similarly, by substituting \refeqn{delq}, \refeqn{we} into \refeqn{dotw}, we obtain
\begin{align}
\delta\dot\omega &= -k_\omega\delta\omega -k_qq_d\times(\xi\times q)\nonumber\\
&=-k_\omega\omega +k_q\hat q_d \hat q\,\xi,\label{eqn:dotdelw}
\end{align}
which is the linearized equation for \refeqn{dotw}. 

Equations \refeqn{dotxi}, \refeqn{dotdelw} can be written in a matrix form as
\begin{align}
\dot x =
\begin{bmatrix}\dot\xi \\ \delta\dot\omega \end{bmatrix}
=\begin{bmatrix} qq^T\hat\omega & I-qq^T\\k_q\hat q_d\hat q & -k_wI\end{bmatrix}
\begin{bmatrix}\xi \\ \delta\omega \end{bmatrix}=Ax,\label{eqn:xdot}
\end{align}
where the state vector of the linearized controlled system is $x=[\xi;\delta\omega]\in\Re^6$. A spherical pendulum has two degrees of freedom, but this linearized equation of motion evolves in $\Re^6$ instead of $\Re^4$. Since $q\cdot\omega=0$ and $q\cdot\xi=0$, we have the following two additional constraints on $\xi,\delta\omega$:
\begin{align}
Cx=
\begin{bmatrix} q^T & 0 \\ -\omega^T\hat q & q^T\end{bmatrix}
\begin{bmatrix}\xi \\ \delta\omega \end{bmatrix}
=\begin{bmatrix} 0 \\ 0 \end{bmatrix}.\label{eqn:con}
\end{align}
Therefore, the state vector $x$ should lie in the null space of the matrix $C\in\Re^{2\times 4}$. However, this is not an extra constraint that should be imposed when solving \refeqn{xdot}. As long as the initial condition $x(0)$ satisfies \refeqn{con}, the structure of \refeqn{dotw}, \refeqn{dotq}, and \refeqn{xdot}, guarantees that the state vector $x(t)$ satisfies \refeqn{con} for all $t$, i.e. $\frac{d}{dt}C(t)x(t) =0$ for all $t\geq 0$ when $C(0)x(0)=0$. This means that the null space of $C$ is a flow-invariant subspace.

\subsection{Equilibrium Solutions}

We choose the desired direction as $q_d=e_3$. The equilibrium solution $(q_d,0)=(e_3,0)$ is referred to as the hanging equilibrium, and the additional equilibrium solution $(-q_d,0)=(-e_3,0)$ is referred to as the inverted equilibrium. We study the eigen-structure of each equilibrium using the linearized equation \refeqn{xdot}. To illustrate the ideas, the controller gains are selected as $k_q=k_\omega=1$. 

\subsubsection{Hanging Equilibrium} 

The eigenvalues $\lambda_i$, and the eigenvectors $v_i$ of the matrix $A$ at the hanging equilibrium $(e_3,0)$ are given by
\begin{gather*}
\lambda_{1,2}=(-1\pm\sqrt{3}i)/2,\;
\lambda_{3,4}=\lambda_{1,2},\;\lambda_5=0,\;\lambda_6=-1,\\
v_{1,2}= e_1 + (-1\pm\sqrt{3}i)e_4/2,\;
v_{3,4}= e_2 + (-1\pm\sqrt{3}i)e_5/2,\\
v_5=e_3,\quad v_6=e_6,
\end{gather*}
where $e_i\in\Re^6$ denotes the unit-vector whose $i$-th element is one, and other elements are zeros. Note that there are repeated eigenvalues, but we obtain six linearly independent eigenvectors, i.e., the geometric multiplicities are equal to the algebraic multiplicities.

The basis of the null space of the matrix $C$, namely $\mathcal{N}(C)$ is $\{e_1,e_2,e_4,e_5\}$. The solution of the linearized equation can be written as $x(t)=\sum_{i=1}^6 c_i \exp(\lambda_it) v_i$ for constants $c_i$ that are determined by the initial condition: $x(0)=\sum_{i=1}^6 c_i v_i$. But, the eigenvectors $v_5,v_6$ do not satisfy the constraint given by \refeqn{con}, since they do not lie in $\mathcal{N}(C)$. Therefore, the constants $c_5,c_6$ are zero for initial conditions that are compatible with \refeqn{con}. We have $\mathrm{Re}[\lambda_i]<0$ for $1\leq i\leq 4$. Therefore, the equilibrium $(q,\omega)=(e_3,0)$ is asymptotically stable.

\subsubsection{Inverted Equilibrium} 

The eigenvalues $\lambda_i$, and the eigenvectors $v_i$ of the matrix $A$ at the inverted equilibrium $(-e_3,0)$ are given by
\begin{gather}
\lambda_{1,2}=-(\sqrt{5}+1)/2,\lambda_{3,4}=(\sqrt{5}-1)/2,\lambda_5=0,\lambda_6=-1,\nonumber\\
v_{1}= e_1 -(\sqrt{5}+1)e_4/2,\, v_2=e_2-(\sqrt{5}+1)e_5/2,\label{eqn:v1v2}\\
v_3=(\sqrt{5}+1)e_1/2 +e_4,\, v_4=(\sqrt{5}+1)e_2/2+e_5,\nonumber\\
v_5=e_3,\, v_6=e_6.\nonumber
\end{gather}
The basis of $\mathcal{N}(C)$ is $\{e_1,e_2,e_4,e_5\}$. Hence, the eigenvectors $v_5,v_6$ do not lie in $\mathcal{N}(C)$. Therefore, the solution can be written as $x(t)=\sum_{i=1}^4 c_i \exp(\lambda_it) v_i$ for constants $c_i$ that are determined by the initial condition.

We have $\mathrm{Re}[\lambda_{1,2}]<0$, and $\mathrm{Re}[\lambda_{3,4}]>0$. Therefore, the inverted equilibrium $(q,\omega)=(-e_3,0)$ is a hyperbolic equilibrium, and in particular, a saddle point.

\subsection{Stable Manifold for the Inverted Equilibrium}\label{sec:SM}

\subsubsection{Stable Manifold}

The saddle point $(-e_3,0)$ has a stable manifold $W^s$, which is defined to be
\begin{align*}
W^s(-e_3,0) &= \{ (q,\omega)\in \T\Sph^2\,|\, \lim_{t\rightarrow\infty} \mathcal{F}^t(q,\omega) = (-e_3,0)\},
\end{align*}
where $\mathcal{F}^t:(q(0),\omega(0))\rightarrow(q(t),\omega(t))$ denotes the flow map along the solution of \refeqn{dotw}, \refeqn{dotq}. The existence of $W^s(-e_3,0)$ has nontrivial effects on the overall dynamics of the controlled system. Trajectories in $W^s(-e_3,0)$ converge to the antipodal point of the desired equilibrium $(e_3,0)$, and it takes a long time period for any trajectory near $W^s(-e_3,0)$ to asymptotically converge to the desired equilibrium $(e_3,0)$. 

According to the stable and unstable manifold theorem~\cite{Kuz98}, a local stable manifold $W^s_{loc}(-e_3,0)$ exists in the neighborhood of $(-e_3,0)$, and it is tangent to the stable eigenspace $E^s(-e_3,0)$ spanned by the eigenvectors $v_1$ and $v_2$ of the stable eigenvalues $\lambda_{1,2}$. The (global) stable manifold can be written as
\begin{align}
W^s(-e_3,0) & = \bigcup_{t>0} \mathcal{F}^{-t} ( W^s_{loc}(-e_3,0)),\label{eqn:Ws}
\end{align}
which states that the stable manifold $W^s$ can be obtained by globalizing the local stable manifold $W^s_{loc}$ by the backward flow map. 

This yields a method to compute $W^s(-e_3,0)$~\cite{KraOsiIJBC05}. We choose a small ball $B_\delta\subset W^s_{loc}(-e_3,0)$ with a radius $\delta$ around $(-e_3,0)$, and we grow the manifold $W^s(-e_3,0)$ by evolving $B_\delta$ under the flow $\mathcal{F}^{-t}$. More explicitly, the stable manifold can be parameterized by $t$ as follows:
\begin{align}
W^s(-e_3,0) =\{ \mathcal{F}^{-t} (B_\delta)\}_{t>0}.\label{eqn:Wc}
\end{align}

We construct a ball in the stable eigenspace of $(-e_3,0)$ with sufficiently small radius $\delta$, i.e. $B_\delta\subset E^s_{loc}(-e_3,0)$. From the stable eigenvectors $v_1,v_2$ at \refeqn{v1v2}, $E^s_{loc}(-e_3,0)$ can be written as
\begin{align}
E^s_{loc}& (-e_3,0) = \{  (q,\omega)\in \T\Sph^2\,|\,
  q=\exp(\alpha_1\hat e_1+\alpha_2\hat e_2)(-e_3),\nonumber\\
& \omega=-\hat q^2(-(\sqrt{5}+1)/2)(\alpha_1 e_1+\alpha_2e_2)\text{ for $\alpha_1,\alpha_2\in\Re$}\},\label{eqn:Esloc}
\end{align}
where $-\hat q^2$ in the expression for $\omega$ corresponds to the orthogonal projection onto the plane normal to $q$, as required due to the constraint $q\cdot\omega=0$. 

We define a distance on $\T\Sph^2$ as follows:
\begin{align}
d_{\T\Sph^2} ((q_1,\omega_1),(q_2,\omega_2)) = \sqrt{\Psi(q_1,q_2)} + \|\omega_1-\omega_2\|.\label{eqn:dis}
\end{align}
For $\delta>0$, the subset $B_\delta$ of $E^s_{loc}(-e_3,0)$ is parameterized by $\theta\in\Sph^1$ as
\begin{align}
B_\delta & = \{  (q,\omega)\in \T\Sph^2\,|\,
  q=\exp(\alpha_1\hat e_1+\alpha_2\hat e_2)(-e_3),\nonumber\\
& \omega=-\hat q^2(-(\sqrt{5}+1)/2)(\alpha_1 e_1+\alpha_2e_2),\text{ where}\nonumber\\
& \text{$\alpha_1=\frac{\delta}{1/\sqrt{2}+(\sqrt{5}+1)/2}\cos\theta,$\;}\nonumber\\
& \text{$\alpha_2=\frac{\delta}{1/\sqrt{2}+(\sqrt{5}+1)/2}\sin\theta$,} \text{ for $\theta\in\Sph^1$} \}.\label{eqn:Bdelta}
\end{align}
The given choice of the constants $\alpha_1,\alpha_2$ guarantees that any point in $B_\delta$ has a distance $\delta$ to $(-e_3,0)$ according to the distance metric \refeqn{dis}.

\subsubsection{Variational Integrators}

The parameterization of the stable manifold $W_s$ in \refeqn{Wc} requires the computation of the backward flow map $\mathcal{F}^{-t}$. However, general purpose numerical integrators may not preserve the structure of the two-sphere or the underlying dynamic characteristics, such as energy dissipation rate, accurately, and they may yield qualitatively incorrect numerical results in simulating a complex trajectory over a long-time period~\cite{HaiLub00}. 

Geometric numerical integration is concerned with developing numerical integrators that preserve geometric features of a system, such as invariants, symmetry, and reversibility. In particular, variational integrators are geometric numerical integrators for Lagrangian or Hamiltonian systems, constructed according to Hamilton's principle. They have desirable computational properties of preserving symplecticity and momentum maps, and they exhibit good energy behavior~\cite{MarWesAN01}. A variational integrator is developed for Lagrangian or Hamiltonian systems evolving on the two-sphere in~\cite{LeeLeoIJNME08}. It preserves both the underlying symplectic properties and the structures of the two-sphere concurrently.

A variational integrator on $\Sph^2$ for the controlled dynamics of a spherical pendulum can be written in a backward-time integration form as follows:
\begin{align}
    q_{{k}} & = -\parenth{h\omega_{k+1} - \frac{h^2}{2ml^2} M_{k+1}}\times q_{k+1}\nonumber\\
    &\quad + \parenth{1-\norm{h\omega_{k+1} - \frac{h^2}{2ml^2} M_{k+1}}^2}^{1/2} q_{k+1},\label{eqn:qk}\\
    \omega_{{k}} & = \omega_{k+1} - \frac{h}{2ml^2} M_k -    \frac{h}{2ml^2} M_{k+1},\label{eqn:wk}
\end{align}
where the constant $h>0$ is time step, the subscript $k$ denotes the value of a variable at the time $t_k=kh$, and $M_k = ml^2(-k_\omega \omega_k - k_q q_d\times q_k)$. For given $(q_{k+1},\omega_{k+1})$, we first compute $M_{k+1}$. Then, $q_k$ is obtained by \refeqn{qk}, followed by $M_k$, and $\omega_k$ is computed by \refeqn{wk}. This yields an explicit, discrete inverse flow map $\mathcal{F}_d^{-h}((q_{k+1},\omega_{k+1}))=(q_{k},\omega_{k})$.

\subsubsection{Visualization}

\setlength{\unitlength}{0.1\columnwidth}
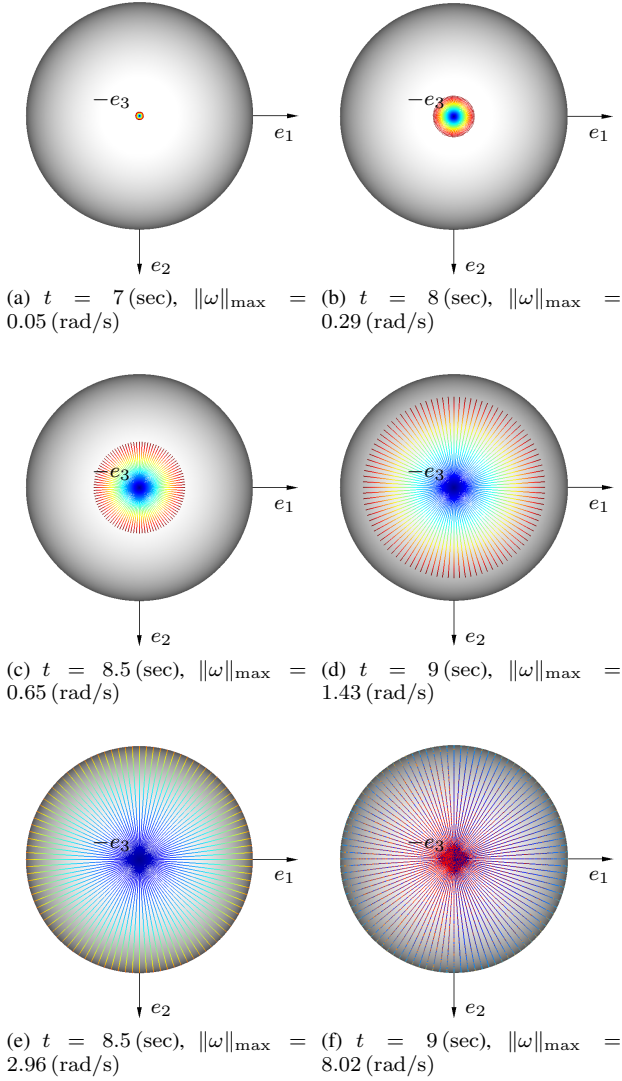
\begin{figure}
\footnotesize\selectfont
\centerline{
\subfigure[$t=7\,(\mathrm{sec}$), $\|\omega\|_{\max}=0.05\,(\mathrm{rad/s})$]{
\begin{picture}(4.5,4.5)(0,0)
\put(0,0){\includegraphics[width=0.45\columnwidth]{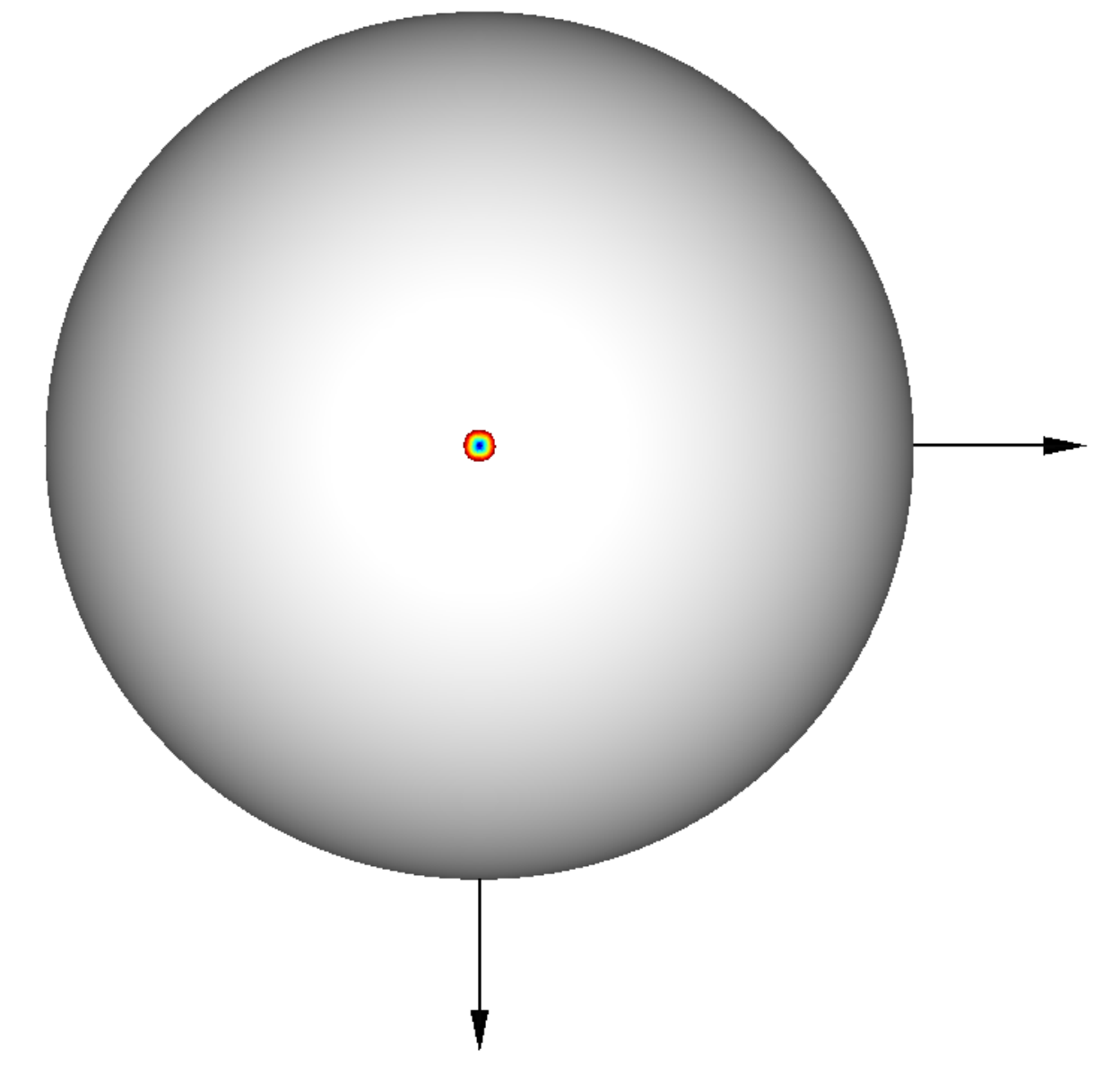}}
\put(4.0,2.2){$e_1$}\put(2.1,0.15){$e_2$}\put(1.2,2.7){$-e_3$}
\end{picture}}
\hspace*{.0cm}
\subfigure[$t=8\,(\mathrm{sec}$), $\|\omega\|_{\max}=0.29\,(\mathrm{rad/s})$]{
\begin{picture}(4.5,4.5)(0,0)
\put(0,0){\includegraphics[width=0.45\columnwidth]{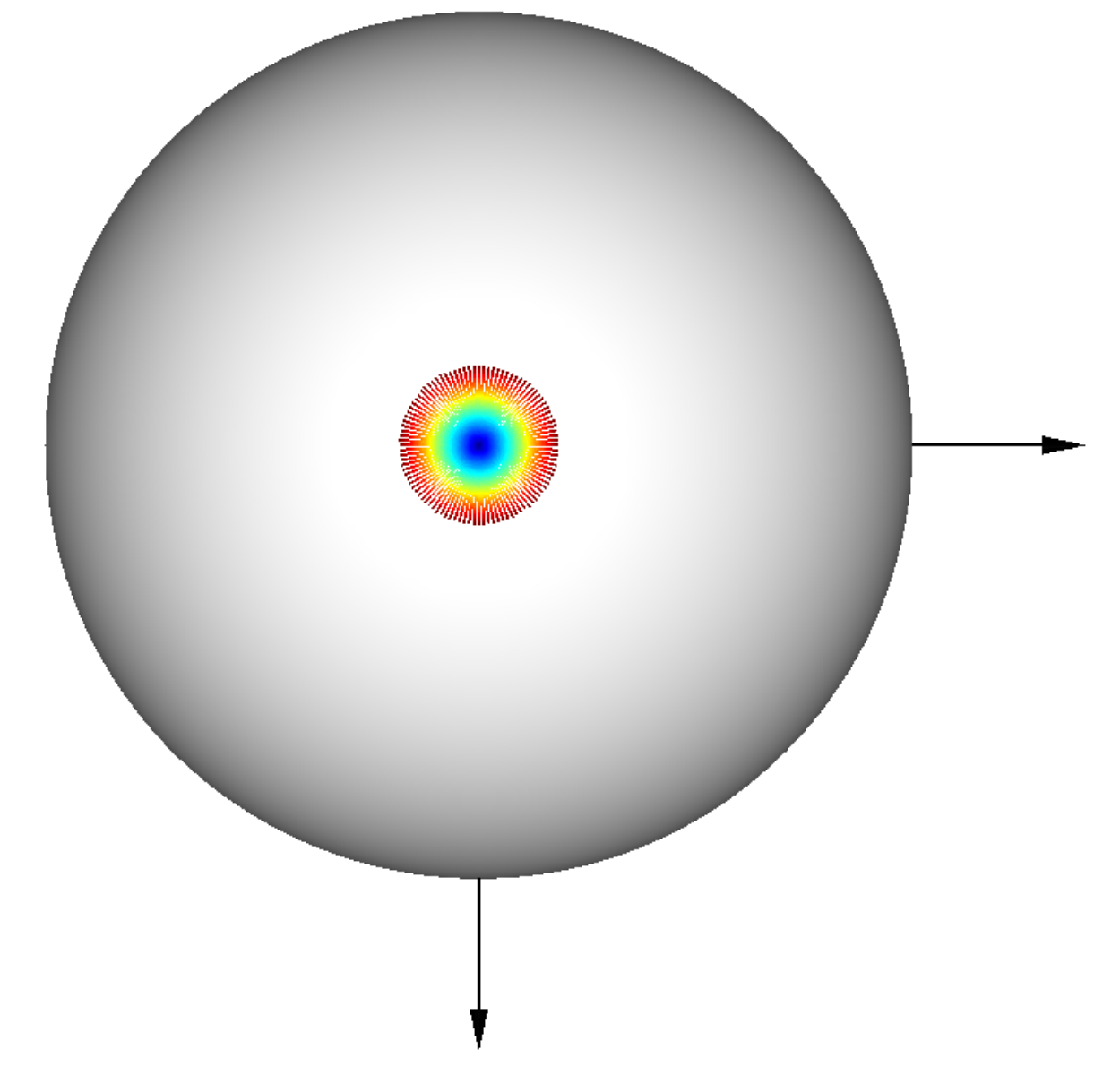}}
\put(4.0,2.2){$e_1$}\put(2.1,0.15){$e_2$}\put(1.2,2.7){$-e_3$}
\end{picture}}
}
\centerline{
\subfigure[$t=8.5\,(\mathrm{sec}$), $\|\omega\|_{\max}=0.65\,(\mathrm{rad/s})$]{
\begin{picture}(4.5,4.5)(0,0)
\put(0,0){\includegraphics[width=0.45\columnwidth]{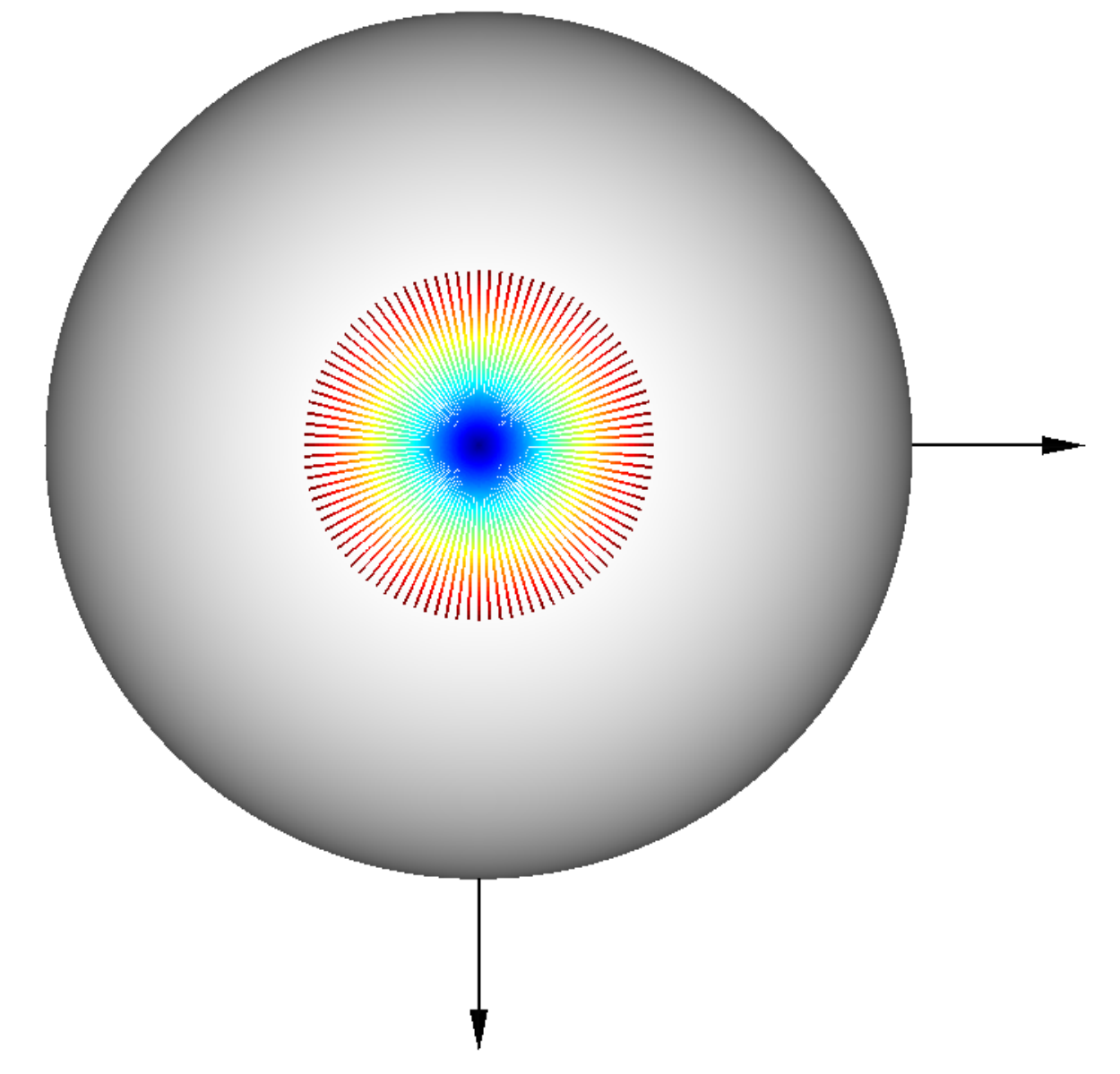}}
\put(4.0,2.2){$e_1$}\put(2.1,0.15){$e_2$}\put(1.2,2.7){$-e_3$}
\end{picture}}
\hspace*{.0cm}
\subfigure[$t=9\,(\mathrm{sec}$), $\|\omega\|_{\max}=1.43\,(\mathrm{rad/s})$]{
\begin{picture}(4.5,4.5)(0,0)
\put(0,0){\includegraphics[width=0.45\columnwidth]{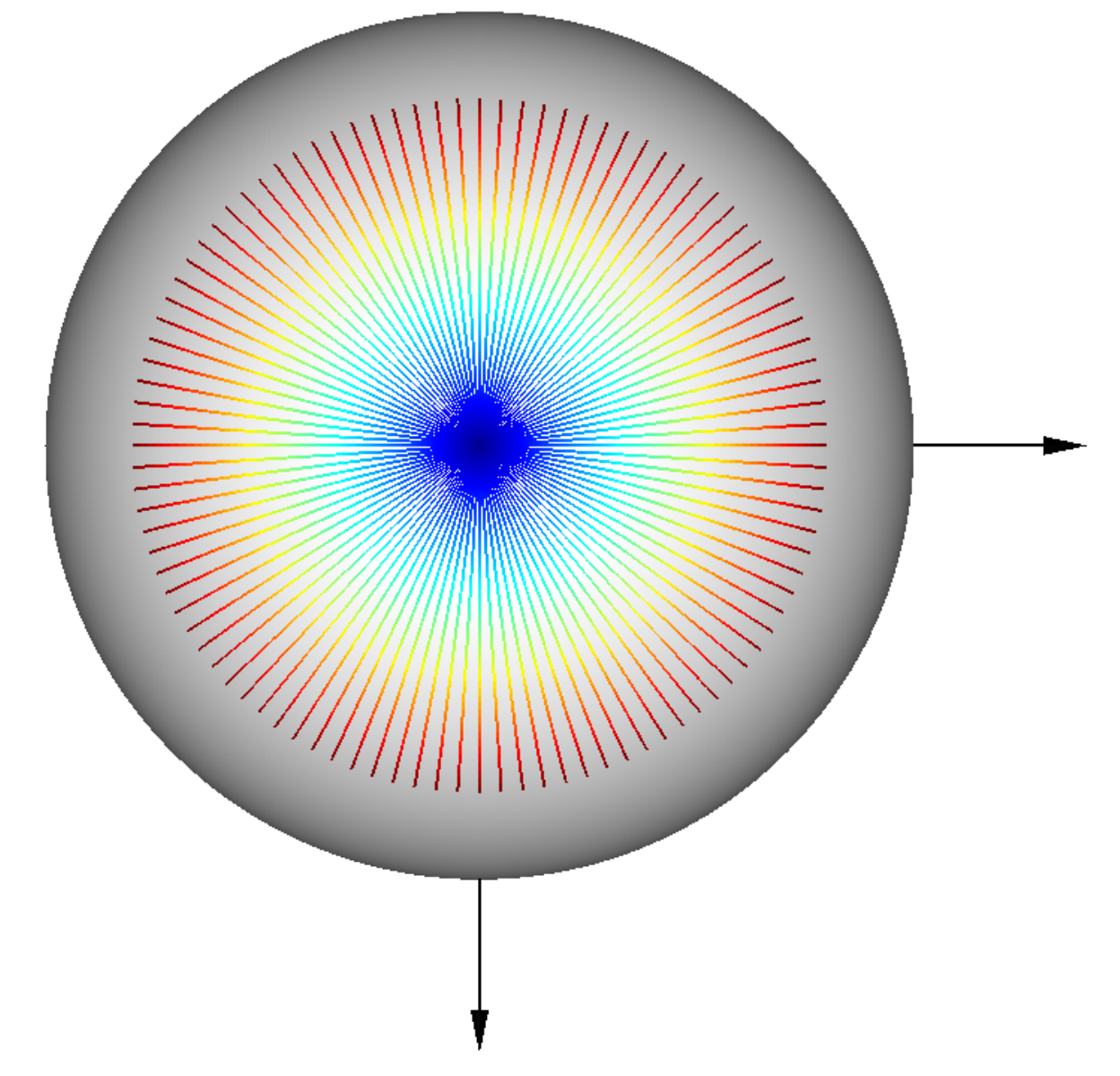}}
\put(4.0,2.2){$e_1$}\put(2.1,0.15){$e_2$}\put(1.2,2.7){$-e_3$}
\end{picture}}
}
\centerline{
\subfigure[$t=8.5\,(\mathrm{sec}$), $\|\omega\|_{\max}=2.96\,(\mathrm{rad/s})$]{
\begin{picture}(4.5,4.5)(0,0)
\put(0,0){\includegraphics[width=0.45\columnwidth]{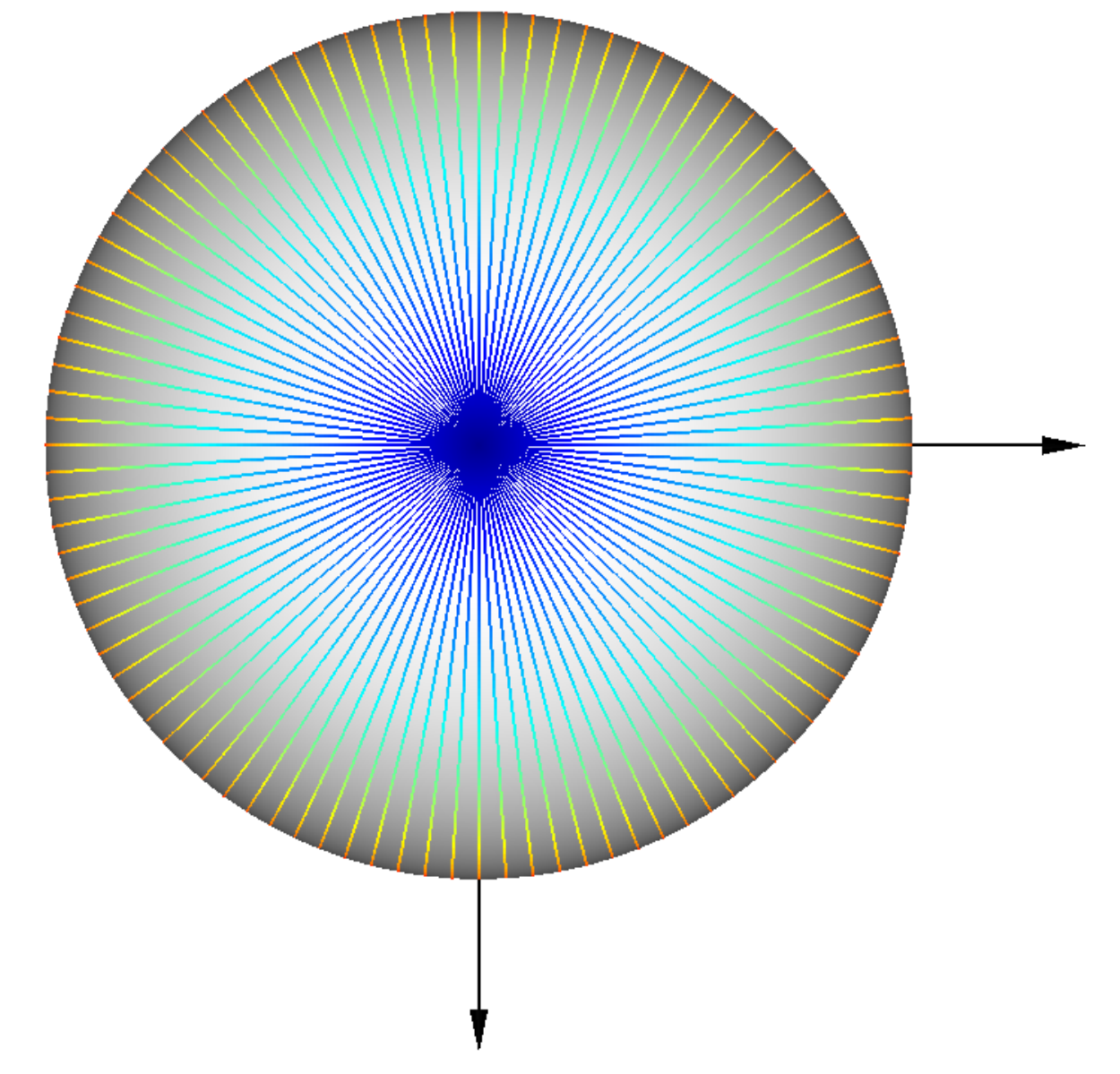}}
\put(4.0,2.2){$e_1$}\put(2.1,0.15){$e_2$}\put(1.2,2.7){$-e_3$}
\end{picture}}
\hspace*{.0cm}
\subfigure[$t=9\,(\mathrm{sec}$), $\|\omega\|_{\max}=8.02\,(\mathrm{rad/s})$]{
\begin{picture}(4.5,4.5)(0,0)
\put(0,0){\includegraphics[width=0.45\columnwidth]{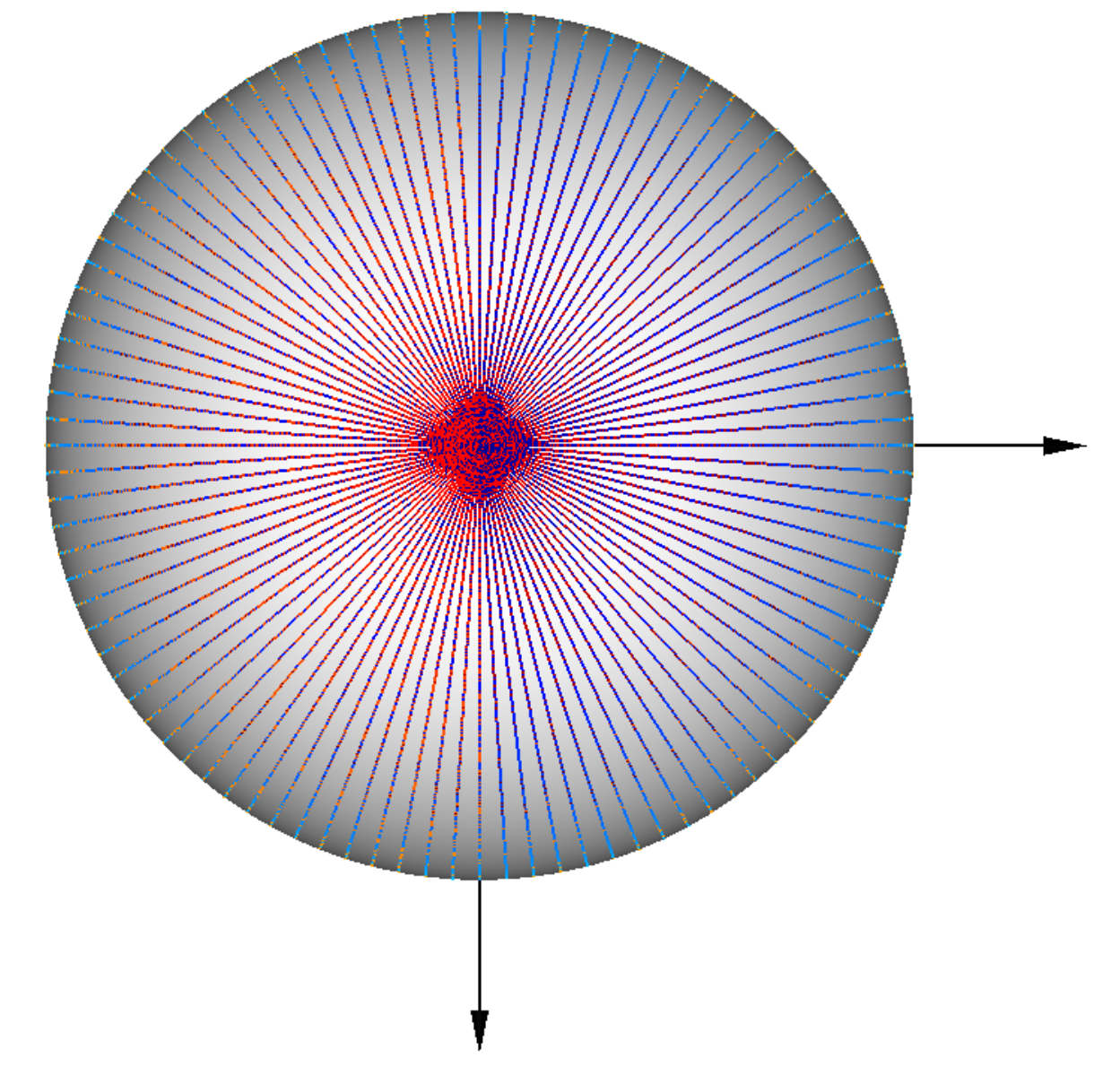}}
\put(4.0,2.2){$e_1$}\put(2.1,0.15){$e_2$}\put(1.2,2.7){$-e_3$}
\end{picture}\label{fig:Ws9}}
}

\caption{Stable manifold to $(q,\omega)=(-e_3,0)$ represented by $\{ \mathcal{F}^{-t} (B_\delta)\}_{t>0}$ for several values of $t$. One hundred points of $B_\delta$ in the stable eigenspace to $(-e_3,0)$ are chosen with $\delta=10^{-6}$, and they are integrated backward in time. Each trajectory is illustrated on a sphere, where the magnitude of angular velocity at each point is denoted by color shading (red: $\|\omega\|_{\max}$, blue: $\|\omega\|_{\min}\simeq 0$).}\label{fig:Ws}
\end{figure}

We choose 100 points on the surface of $B_\delta$ with $\delta=10^{-6}$, and each point is integrated backward using \refeqn{qk}, \refeqn{wk} with timestep $h=0.002$. The resulting trajectories are illustrated in \reffig{Ws} for several values of $t$. Each colored curve on the sphere represents a trajectory on $\T\Sph^2$, since at any point $q$ on the curve, the direction of $\dot q=\omega\times q$ is tangent to the curve at $q$, and the magnitude of $\dot q$ is indirectly represented by color shading. 

We observe the following characteristics of the stable manifold $W_s(-e_3,0)$ of the inverted equilibrium:
\begin{itemize}
\item The boundary of the stable manifold $W_s(-e_3,0)\subset \T\Sph^2$ parameterized by $t$ is circular when projected onto $\Sph^2$. 
\item Each trajectory in $W_s(-e_3,0)$ is on a great circle, when projected onto $\Sph^2$. According to the closed loop dynamics \refeqn{dotw}, and the given initial condition at the surface of $B_\delta$, the direction of $\dot\omega$ is always parallel to $\omega$. Therefore, the direction of $\omega$ is fixed, and the resulting trajectory of $q$ is on a great circle. This also corresponds to the fact that the eigenvalue $\lambda_1$ for the first mode representing the rotations about the first axis is equal to the eigenvalue $\lambda_2$ for the second mode representing the rotations about the second axis at \refeqn{v1v2}, i.e. the convergence rates of these two rotations are identical.
\item The angular velocity decreases to zero as the direction of the pendulum $q$ converges to $-e_3$.
\item The stable manifold $W_s(-e_3,0)$ may cover $\Sph^2$ multiple times if $t$ is sufficiently large, as illustrated at \reffig{Ws9}. Therefore, at any point $q\in\Sph^2$, we can choose $\omega$ such that $(q,\omega)$ lies in the stable manifold $W^s(-e_3,0)$ (the corresponding value of $\omega$ is not unique, since if it is sufficiently large, $q$ can traverse the sphere several times before converging to $-e_3$). This is similar to \textit{kicking} a damped spherical pendulum carefully such that it converges to the inverted equilibrium. 
\end{itemize}

\section{3D Pendulum}

A 3D pendulum is a rigid body supported by a frictionless pivot acting under a gravitational potential. This is a generalization of a planar pendulum or a spherical pendulum, as it has three rotational degrees of freedom. It has been shown that a 3D pendulum may exhibit irregular maneuvers~\cite{ChaLeeJNS11}.

We choose a reference frame, and a body-fixed frame. The origin of the body-fixed frame is located at the pivot point. The attitude of a 3D pendulum is the orientation of the body-fixed frame with respect to the reference frame, and it is described by a rotation matrix representing the linear transformation from the body-fixed frame to the reference frame. The configuration manifold of a 3D pendulum is the special orthogonal group, $\SO=\{R\in\Re^{3\times 3}\,|\, R^T R=I,\mathrm{det}[R]=1\}$.

The equations of motion for a 3D pendulum are given by
\begin{gather*}
J\dot\Omega + \Omega\times J\Omega = mg \rho\times R^T e_3 + u,\label{eqn:Wdot2}\\
\dot R = R\hat\Omega,\label{eqn:Rdot}
\end{gather*}
where the matrix $J\in\Re^{3\times 3}$ is the inertia matrix of the pendulum about the pivot, and $\rho\in\Re^3$ is the vector from the pivot to the center of mass of the pendulum represented in the body-fixed frame. The control moment at the pivot is denoted by $u\in\Re^3$.

\subsection{Control System} 

Several control systems have been developed on $\SO$~\cite{BulLew05,ChaMcCITAC09,LeePACC11}. Here, we summarize a control system to stabilize a 3D pendulum to a fixed desired attitude $R_d\in\SO$. Consider an attitude error function given by
\begin{align*}
\Psi(R,R_d)=\frac{1}{2}\tr{(I-R_d^TR)G},
\end{align*}
for a diagonal matrix $G=\mathrm{diag}[g_1,g_2,g_3]\in\Re^{3\times 3}$ with $g_1,g_2,g_3>0$. The derivative of this attitude error function with respect to $R$ along the direction of $\delta R= R\hat\eta$ for $\eta\in\Re^3$ is given by
\begin{align*}
\D_R & \Psi(R,R_d)\cdot\delta R = -\frac{1}{2}\tr{R_d^TR\hat\eta G} \\
& = \frac{1}{2}(GR_d^TR -R^TR_d G)^\vee \cdot\eta\equiv e_R\cdot \eta,
\end{align*}
where we use the property that $\mathrm{tr}[\hat x A]=-x\cdot(A-A^T)^\vee$ for any $x\in\Re^3,A\in\Re^{3\times 3}$. The \textit{vee map}, $\vee:\so\rightarrow\Re^3$, denotes the inverse of the hat map. An attitude error vector is defined as $e_R = \frac{1}{2}(GR_d^TR -R^TR_d G)\in\Re^3$. For positive constants $k_\Omega,k_R$, we choose the following control input:
\begin{align*}
u = -k_R e_R -k_\Omega \Omega -mg \rho\times R^T e_3.
\end{align*}
The corresponding closed loop dynamics are given by
\begin{gather}
J\dot\Omega =- \Omega\times J\Omega  -k_R e_R -k_\Omega \Omega,\label{eqn:Wdot2}\\
\dot R = R\hat\Omega.\label{eqn:Rdot}
\end{gather}
This system has four equilibria: in addition to the desired equilibrium $(R_d,0)$, there exist three other equilibria at $(R_d\exp (\pi\hat e_i,0),0)$ for $i\in\{1,2,3\}$, which correspond to the rotation of the desired attitude by $180^\circ$ about each body-fixed axis.

The existence of additional, undesirable equilibria is due to the nonlinear topological structure of $\SO$, and it cannot be avoided by constructing a different control system (as long as it is continuous). It has been shown that it is not possible to design a continuous feedback control stabilizing an attitude globally on \SO~\cite{BhaBerSCL00,KodPICDC98}.

The stability of the desired equilibrium can be studied by using the following Lyapunov function,
\begin{align*}
\mathcal{V} =\frac{1}{2}\Omega\cdot J\Omega + k_R \Psi(R,R_d).
\end{align*}
In this paper, we analyze the stability of each equilibrium by linearizing the closed loop dynamics to study the equilibrium structures more explicitly.

\subsection{Linearization} 

A variation in $\SO$ can be expressed as ~\cite{LeeLeoPICCA05}:
\begin{align}
R^\epsilon=R\exp(\epsilon\hat\eta),\quad \Omega^\epsilon=\Omega +\epsilon\delta\Omega,\label{eqn:delRdelW}
\end{align}
for $\eta,\delta\Omega\in\Re^3$. The corresponding infinitesimal variation of $R$ is given by $\delta R = R\hat\eta$. Substituting this into \refeqn{Rdot},
\begin{align*}
R\hat\Omega\hat\eta + R\hat{\dot\eta}= R\hat\eta\hat\Omega + R\delta\hat\Omega.
\end{align*}
Using the property $\hat x \hat y -\hat y\hat x=\widehat{x\times y}$ for any $x,y,\in\Re^3$, this can be rewritten as
\begin{align}
\dot\eta = \delta\Omega -\hat\Omega\eta.\label{eqn:etadot}
\end{align}
Similarly, by substituting \refeqn{delRdelW} into \refeqn{Wdot2}, we obtain
\begin{align}
J\delta\dot\Omega & = -\delta\Omega\times J\Omega -\Omega\times J\delta\Omega\nonumber\\
&\quad -\frac{1}{2}k_R (GR_d^T R\hat\eta +\hat\eta R^T R_d G) -k_\Omega\delta\Omega,\nonumber\\
& = (\widehat{J\Omega}-\hat\Omega J -k_\Omega I) \delta\Omega -\frac{1}{2}k_RH\eta,\label{eqn:delWdot}
\end{align}
where $H=\mathrm{tr}[R^T R_d G]I-R^T R_d G\in\Re^{3\times 3}$, and we used the property, $\hat x A + A^T\hat x = \mathrm{tr}[A]I-A$ for any $x\in\Re^3, A\in\Re^{3\times 3}$. Equations \refeqn{etadot},\refeqn{delWdot} can be written in matrix form as
\begin{align}
\dot x & = \begin{bmatrix}\dot\eta \\ \delta\dot\Omega \end{bmatrix}
=\begin{bmatrix}
-\hat\Omega & I\\
-\frac{1}{2}k_R J^{-1}H & J^{-1}(\widehat{J\Omega}-\hat\Omega J -k_\Omega I)
\end{bmatrix}
\begin{bmatrix}\eta \\ \delta\Omega \end{bmatrix}\nonumber\\
& = Ax.\label{eqn:xdotSO3}
\end{align}
This corresponds to the linearized equation of motion of \refeqn{Wdot2}, \refeqn{Rdot}.

\subsection{Equilibrium Solutions}

We choose the desired attitude as $R_d=I$. In addition to the desired equilibrium $(I,0)$, there are three additional equilibria, namely $(\exp(\pi\hat e_1),0)$, $(\exp(\pi\hat e_2),0)$, $(\exp(\pi\hat e_3),0)$. We study the eigen-structure of each equilibrium using the linearized equation \refeqn{xdotSO3}. We assume that
\begin{align*}
J=\mathrm{diag}[3,2,1]\,\mathrm{kgm^2},\; G=\mathrm{diag}[0.9,1,1.1],\; k_R=k_\Omega=1.
\end{align*}

\subsubsection{Equilibrium $(I,0)$} 
The eigenvalues of the matrix $A$ at the desired equilibrium $(I,0)$ are given by
\begin{gather*}
\lambda_{1,2}=-0.1667\pm0.5676i,\\
\lambda_{3,4}=-0.25\pm 0.6614i,\\
\lambda_{5,6}=-0.5\pm 0.8367i.
\end{gather*}
This equilibrium is an asymptotically stable focus.

\subsubsection{Equilibrium $(\exp(\pi\hat e_1),0)$} 
At this equilibrium, the eigenvalues and the eigenvectors of $A$ are given by
\begin{gather}
\lambda_1=-0.7813,\quad v_1= e_1-0.7813 e_4,\nonumber\\
\lambda_2=-0.5854,\quad v_2=e_2-0.5854e_5,\nonumber\\
\lambda_3=-1.0477,\quad v_3=e_3-1.0477e_6,\label{eqn:v3_SO31}\\
\lambda_4=0.4480,\quad v_4=e_1+0.4480e_4,\nonumber\\
\lambda_5=0.0854,\quad v_5=e_2+0.0854e_5,\nonumber\\
\lambda_6=0.0477,\quad v_6=e_3+0.0477e_6.\nonumber
\end{gather}
Therefore, this equilibrium is a saddle point, where three modes are stable, and three modes are unstable.

\subsubsection{Equilibrium $(\exp(\pi\hat e_2),0)$} 
At this equilibrium, the eigenvalues and the eigenvectors of $A$ are given by
\begin{gather}
\lambda_1=-0.3775,\quad v_1= e_1-0.3775 e_4,\nonumber\\
\lambda_2=-1,\quad v_2=e_2-e_5,\label{eqn:v2_SO32}\\
\lambda_3=-0.9472,\quad v_3=e_3-0.9472e_6,\nonumber\\
\lambda_4=-0.0528,\quad v_4=e_3-0.0528e_6,\nonumber\\
\lambda_5=0.0442,\quad v_5=e_1+0.0442e_4,\nonumber\\
\lambda_6=0.5,\quad v_6=e_2+5e_5.\nonumber
\end{gather}
Therefore, this equilibrium is a saddle point, where four modes are stable, and two modes are unstable.

\subsubsection{Equilibrium $(\exp(\pi\hat e_3),0)$} 
At this equilibrium, the eigenvalues and the eigenvectors of $A$ are given by
\begin{gather}
\lambda_1=-0.0613,\quad v_1= e_1-0.0613 e_4,\nonumber\\
\lambda_2=-0.2721,\quad v_2= e_1-0.2721 e_4,\nonumber\\
\lambda_3=-0.1382,\quad v_3= e_2-0.1382 e_5,\nonumber\\
\lambda_4=-0.3618,\quad v_4= e_2-0.3618 e_5,\nonumber\\
\lambda_5=-1.5954,\quad v_5= e_3-1.5954 e_6,\label{eqn:v5_SO33}\\
\lambda_6= 0.5954,\quad v_6= e_2+0.5954 e_6.\nonumber
\end{gather}
Therefore, this equilibrium is a saddle point, where five modes are stable, and one mode is unstable.

\subsection{Stable Manifolds for the Saddle Points} 

The eigen-structure analysis shows that there exist multi-dimensional stable manifolds for each saddle point. They have zero measure as the dimension of stable manifold is less than the dimension of $\T\SO$. But, the existence of these stable manifolds may have nontrivial effects on the attitude dynamics.

We numerically characterize these stable manifolds using backward time integration, as discussed in Section \ref{sec:SM}.

The stable eigenspace for each saddle point can be written as
\begin{align*}
&E^s_{loc} (\exp(\pi\hat e_1),0) = \{ (R,\Omega)\in \T\SO\,|\,\nonumber\\
&\; R=\exp(\pi\hat e_1)\exp(\alpha_1\hat e_1+\alpha_2\hat e_2+\alpha_3\hat e_3),\\
&\; \Omega=-0.7813\alpha_1e_1-0.5854\alpha_2e_2-1.0477\alpha_3e_3\text{ for $\alpha_i\in\Re$}\},\\
&E^s_{loc} (\exp(\pi\hat e_2),0) = \{ (R,\Omega)\in \T\SO\,|\,\nonumber\\
&\; R=\exp(\pi\hat e_2)\exp(\alpha_1\hat e_1+\alpha_2\hat e_2+(\alpha_3+\alpha_4)\hat e_3),\\
&\; \Omega=-0.37\alpha_1e_1-\alpha_2e_2-(0.94\alpha_3+0.05\alpha_4)e_3\text{ for $\alpha_i\in\Re$}\},\\
&E^s_{loc} (\exp(\pi\hat e_3),0) = \{ (R,\Omega)\in \T\SO\,|\,\nonumber\\
&\; R=\exp(\pi\hat e_3)\exp((\alpha_1+\alpha_2)\hat e_1+(\alpha_3+\alpha_4)\hat e_2+\alpha_5\hat e_3),\\
&\; \Omega=-(0.06\alpha_1+0.27\alpha_2)e_1-(0.13\alpha_3+0.36\alpha_4)e_2\\
&\;\quad -1.59\alpha_5e_3\text{ for $\alpha_i\in\Re$}\},
\end{align*}

We define a distance on $\T\SO$ as follows:
\begin{align*}
d_{\T\SO} ((R_1,\Omega_1),(R_2,\Omega_2)) = \sqrt{\Psi(R_1,R_2)} + \|\Omega_1-\Omega_2\|.
\end{align*}

A variational integrator for the attitude dynamics of a rigid body on $\SO$ is developed in~\cite{LeeLeoPICCA05,LeeLeoCMDA07}. It can be rewritten in a backward integration form as follows:
\begin{gather}
h (\Pi_{k+1}-\frac{h}{2}M_{k+1})^\wedge = J_dF_k - F_k^TJ_d,\label{eqn:Fk}\\
R_k = R_{k+1}F_k^T,\label{eqn:Rk}\\
\Pi_k = F_k\Pi_{k+1}-\frac{h}{2}F_kM_{k+1}-\frac{h}{2} M_k,\label{eqn:Pik}
\end{gather}          
where $M_{k}=u_k+mg\rho\times R^T e_3\in\Re^3$ is the external moment, $\Pi_k=J\Omega_k\in\Re^3$ is the angular momentum. The matrix $J_d\in\Re^{3\times 3}$ denotes a non-standard inertia matrix given by $J_d = \frac{1}{2}\mathrm{tr}[J]I-J$, and the rotation matrix $F_k\in\SO$ represent the relative attitude update between two integration time steps. For given $(R_{k+1},\Pi_{k+1})$, we first compute $M_{k+1}$, and solve \refeqn{Fk} for $F_k$. Then, $R_k$ is obtained by \refeqn{Rk}, and $\Pi_k$ is computed by \refeqn{Pik}. This yields a discrete  inverse flow map, $\mathcal{F}^{-h}_d(R_{k+1},\Pi_{k+1})\rightarrow(R_{k},\Pi_{k})$.

\setlength{\unitlength}{0.1\columnwidth}
\begin{figure*}
\footnotesize\selectfont
\centerline{
\subfigure[$t=11\,(\mathrm{sec}$), $\|\Omega\|_{\max}=0.06\,(\mathrm{rad/s})$]{
\begin{picture}(4.5,4.5)(0,0)
\put(0,0){\includegraphics[width=0.45\columnwidth]{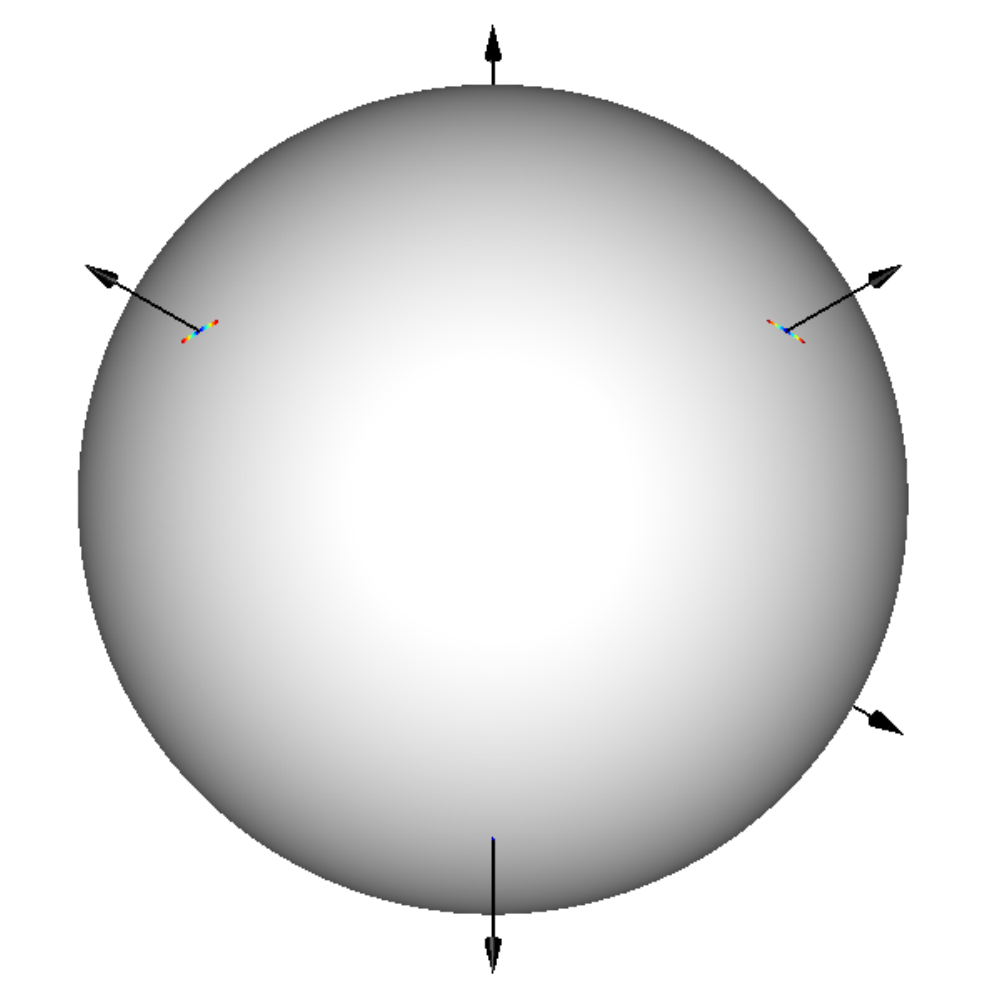}}
\put(3.9,0.9){$e_2$}\put(3.9,3.4){$e_1$}\put(2.3,4.25){$e_3$}
\put(0.1,3.4){$-e_2$}\put(2.3,0.1){$-e_3$}
\end{picture}}
\hspace*{0.2cm}
\subfigure[$t=12\,(\mathrm{sec}$), $\|\Omega\|_{\max}=0.17\,(\mathrm{rad/s})$]{
\begin{picture}(4.5,4.5)(0,0)
\put(0,0){\includegraphics[width=0.45\columnwidth]{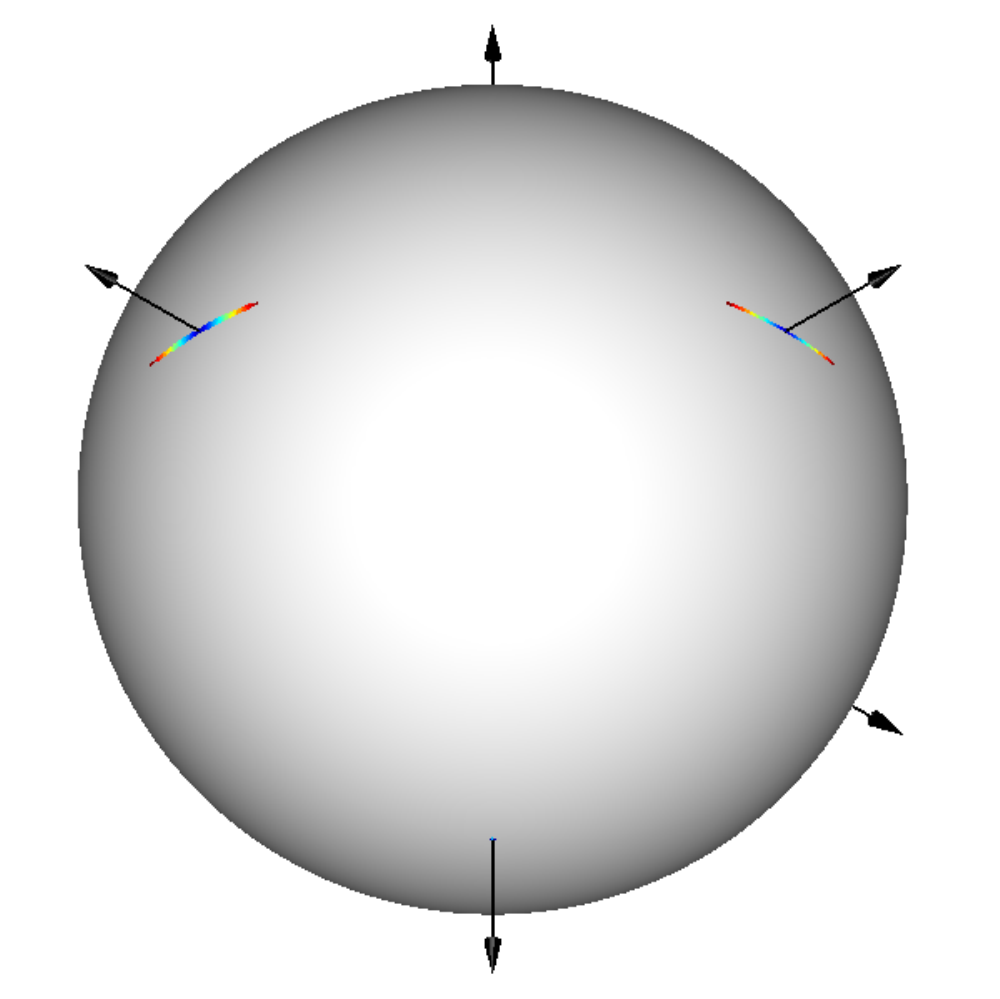}}
\put(3.9,0.9){$e_2$}\put(3.9,3.4){$e_1$}\put(2.3,4.25){$e_3$}
\put(0.1,3.4){$-e_2$}\put(2.3,0.1){$-e_3$}
\end{picture}}
\hspace*{0.2cm}
\subfigure[$t=13\,(\mathrm{sec}$), $\|\Omega\|_{\max}=0.50\,(\mathrm{rad/s})$]{
\begin{picture}(4.5,4.5)(0,0)
\put(0,0){\includegraphics[width=0.45\columnwidth]{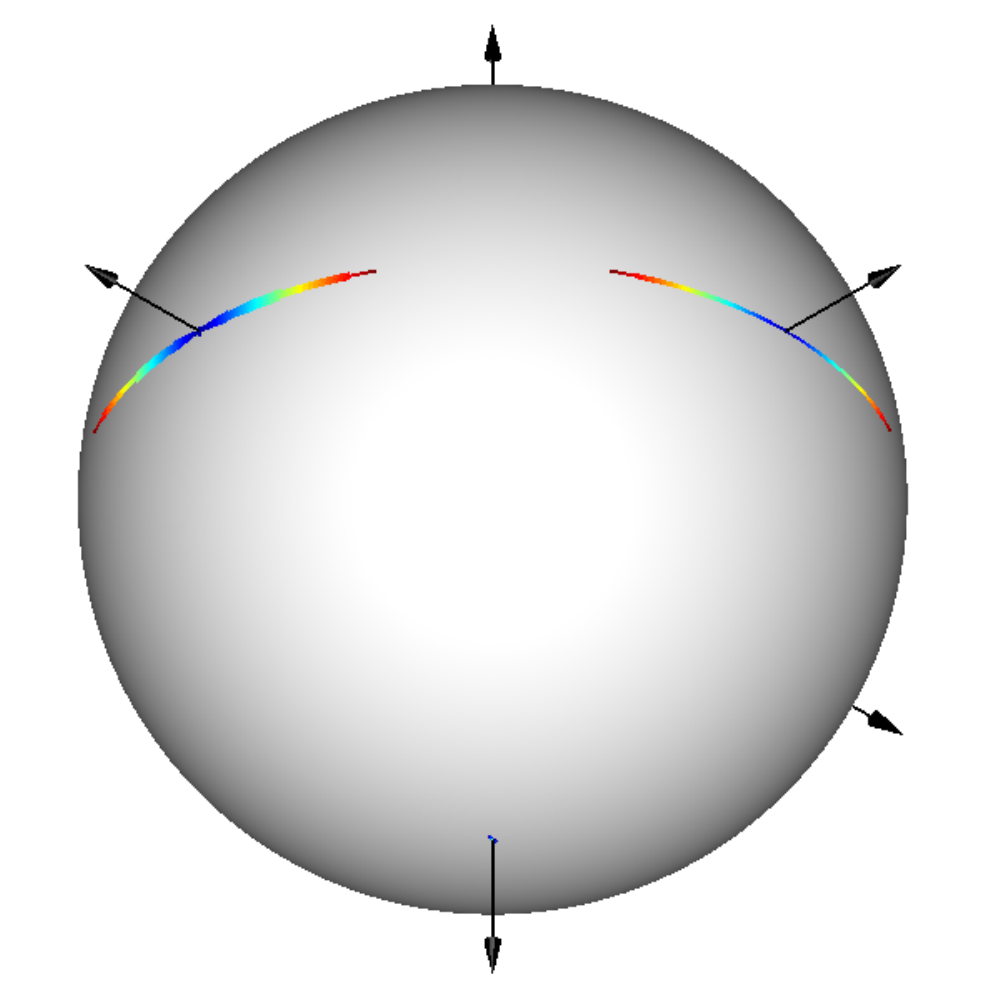}}
\put(3.9,0.9){$e_2$}\put(3.9,3.4){$e_1$}\put(2.3,4.25){$e_3$}
\put(0.1,3.4){$-e_2$}\put(2.3,0.1){$-e_3$}
\end{picture}}
\hspace*{0.2cm}
\subfigure[$t=14\,(\mathrm{sec}$), $\|\Omega\|_{\max}=1.42\,(\mathrm{rad/s})$]{
\begin{picture}(4.5,4.5)(0,0)
\put(0,0){\includegraphics[width=0.45\columnwidth]{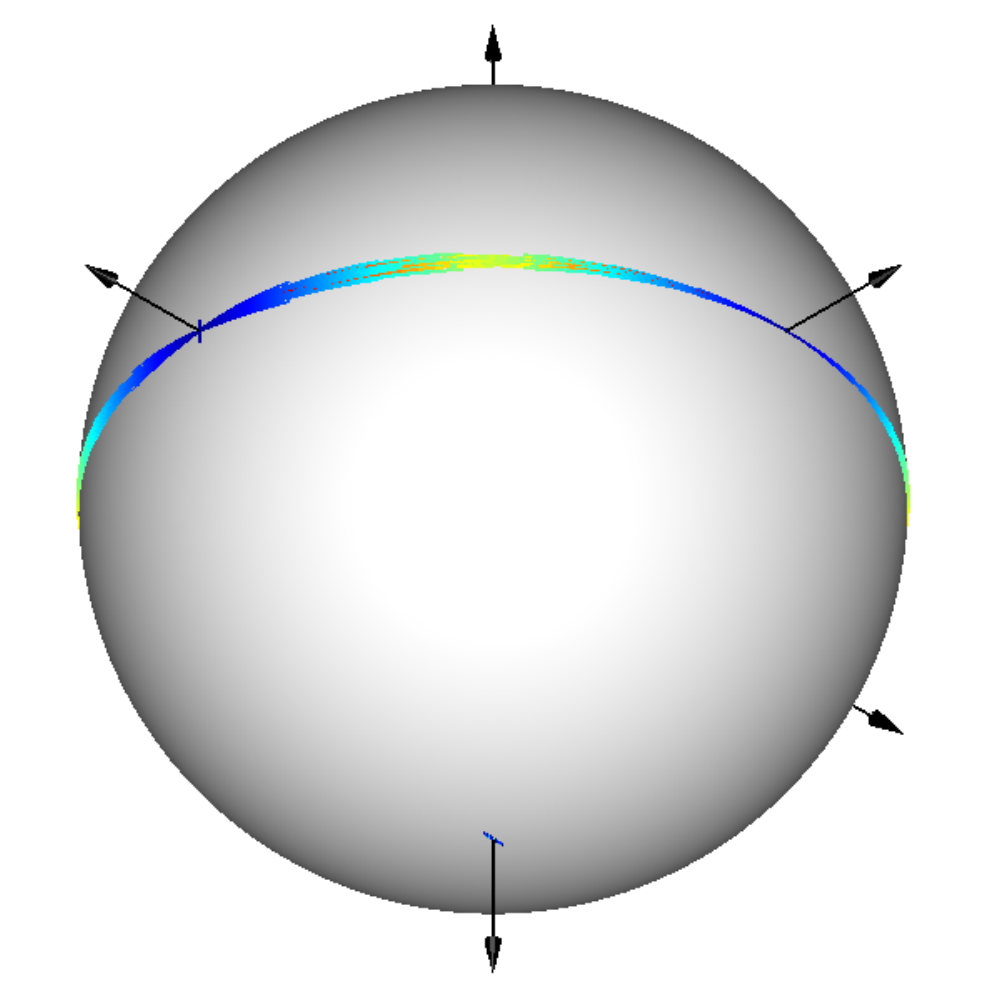}}
\put(3.9,0.9){$e_2$}\put(3.9,3.4){$e_1$}\put(2.3,4.25){$e_3$}
\put(0.1,3.4){$-e_2$}\put(2.3,0.1){$-e_3$}
\end{picture}}
}
\centerline{
\subfigure[$t=15\,(\mathrm{sec}$), $\|\Omega\|_{\max}=3.93\,(\mathrm{rad/s})$]{
\begin{picture}(4.5,4.5)(0,0)
\put(0,0){\includegraphics[width=0.45\columnwidth]{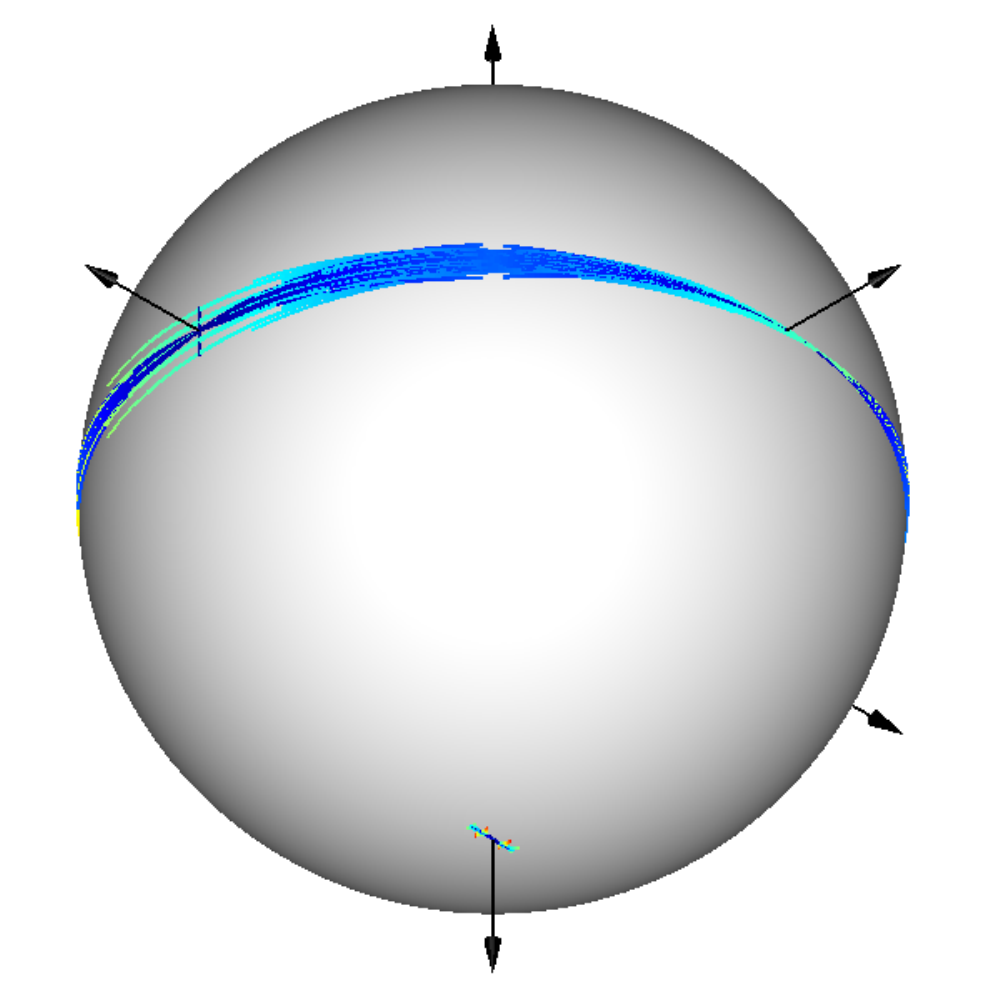}}
\put(3.9,0.9){$e_2$}\put(3.9,3.4){$e_1$}\put(2.3,4.25){$e_3$}
\put(0.1,3.4){$-e_2$}\put(2.3,0.1){$-e_3$}
\end{picture}}
\hspace*{0.2cm}
\subfigure[$t=16\,(\mathrm{sec}$), $\|\Omega\|_{\max}=10.67\,(\mathrm{rad/s})$]{
\begin{picture}(4.5,4.5)(0,0)
\put(0,0){\includegraphics[width=0.45\columnwidth]{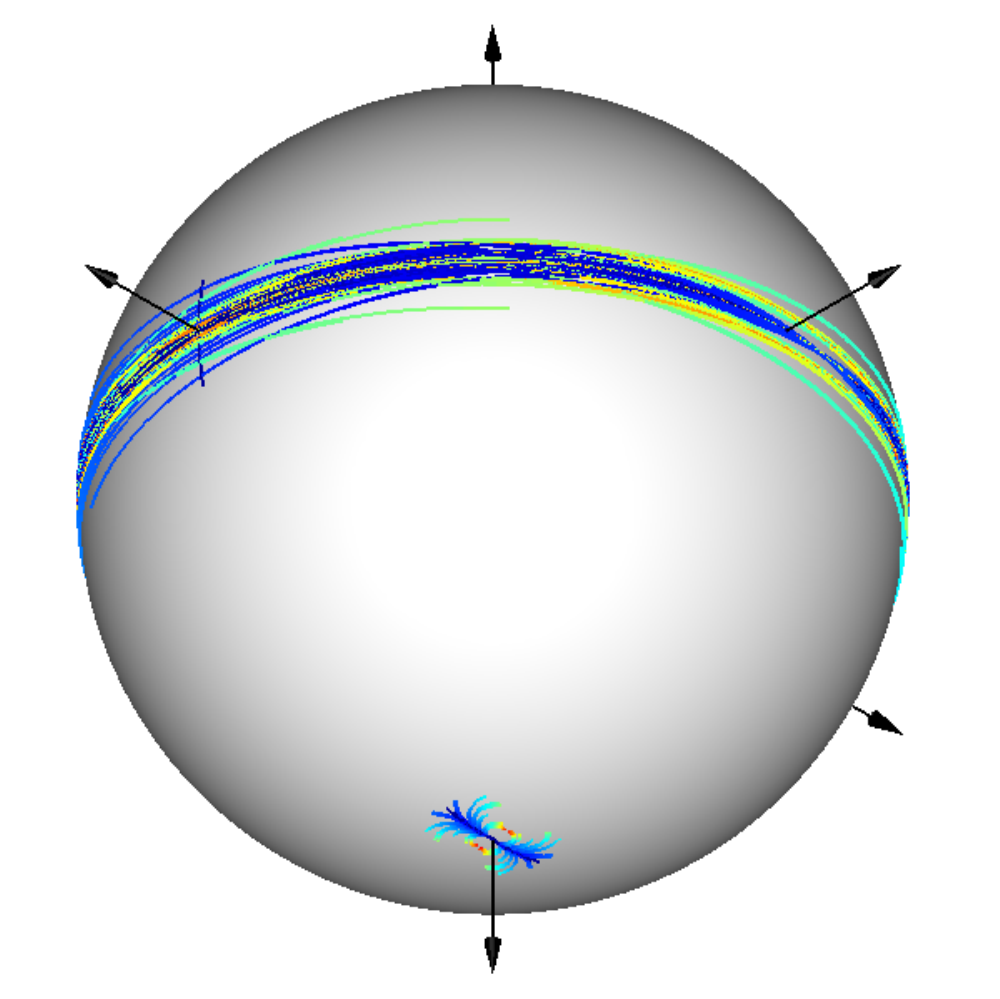}}
\put(3.9,0.9){$e_2$}\put(3.9,3.4){$e_1$}\put(2.3,4.25){$e_3$}
\put(0.1,3.4){$-e_2$}\put(2.3,0.1){$-e_3$}
\end{picture}}
\hspace*{0.2cm}
\subfigure[$t=17\,(\mathrm{sec}$), $\|\Omega\|_{\max}=29.00\,(\mathrm{rad/s})$]{
\begin{picture}(4.5,4.5)(0,0)
\put(0,0){\includegraphics[width=0.45\columnwidth]{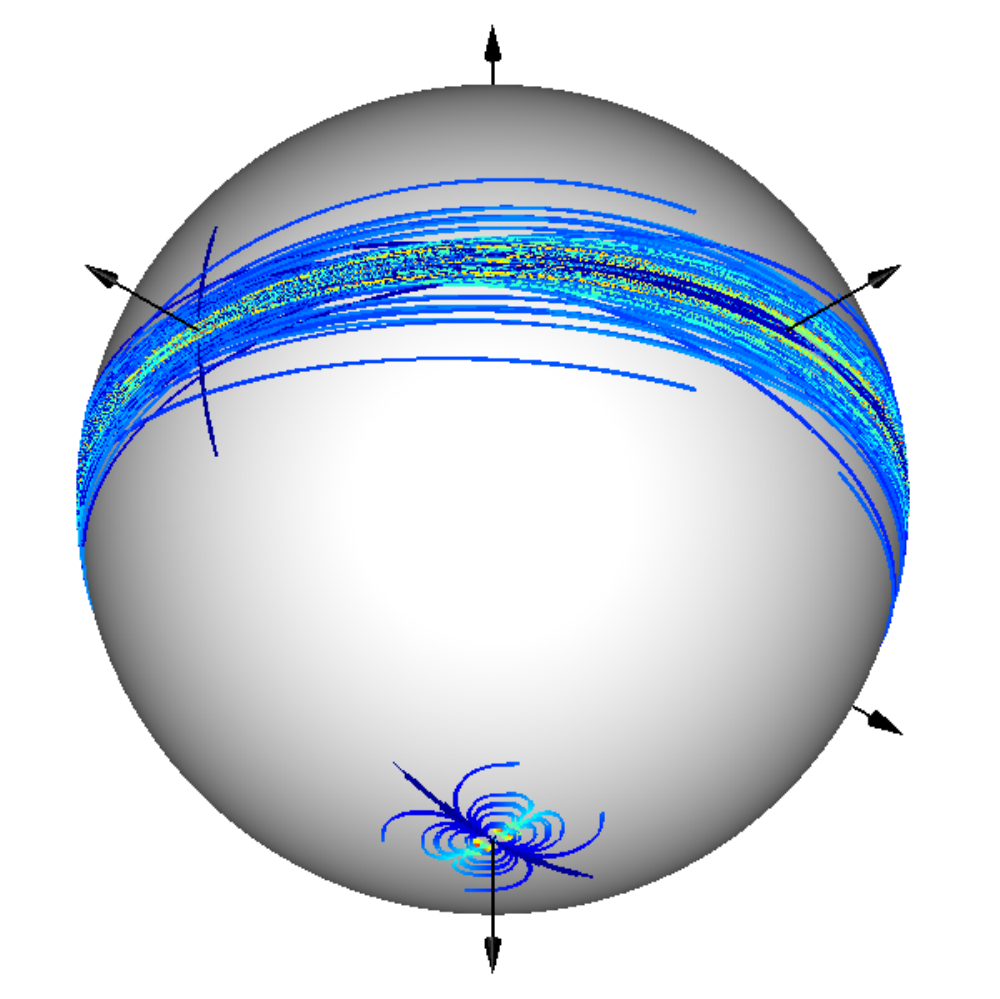}}
\put(3.9,0.9){$e_2$}\put(3.9,3.4){$e_1$}\put(2.3,4.25){$e_3$}
\put(0.1,3.4){$-e_2$}\put(2.3,0.1){$-e_3$}
\end{picture}}
\hspace*{0.2cm}
\subfigure[$t=18\,(\mathrm{sec}$), $\|\Omega\|_{\max}=78.84\,(\mathrm{rad/s})$]{
\begin{picture}(4.5,4.5)(0,0)
\put(0,0){\includegraphics[width=0.45\columnwidth]{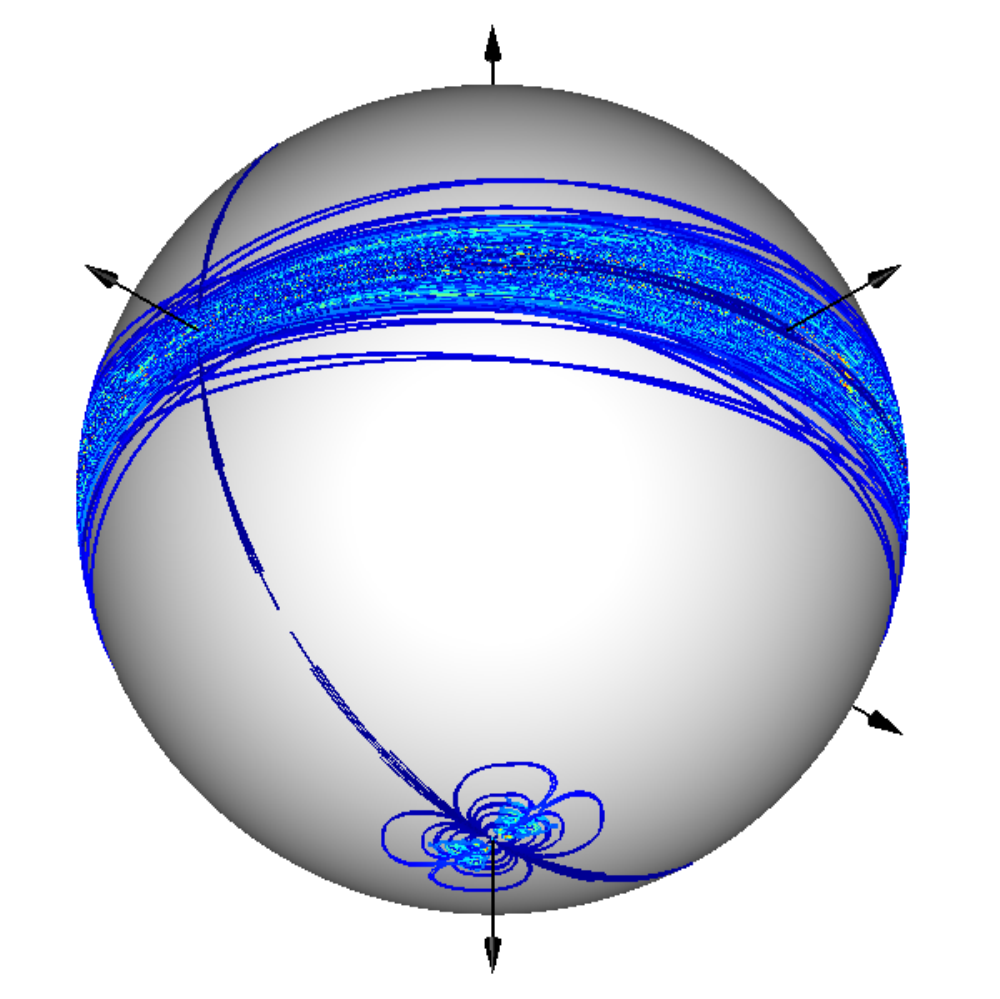}}
\put(3.9,0.9){$e_2$}\put(3.9,3.4){$e_1$}\put(2.3,4.25){$e_3$}
\put(0.1,3.4){$-e_2$}\put(2.3,0.1){$-e_3$}
\end{picture}}
}
\caption{Stable manifold to $(\exp(\pi\hat e_1),0)=([e_1,-e_2,-e_3],0)$ represented by $\{ \mathcal{F}^{-t} (B_\delta)\}_{t>0}$ with $\delta=10^{-6}$ for several values of $t$.}\label{fig:SM1}
\end{figure*}
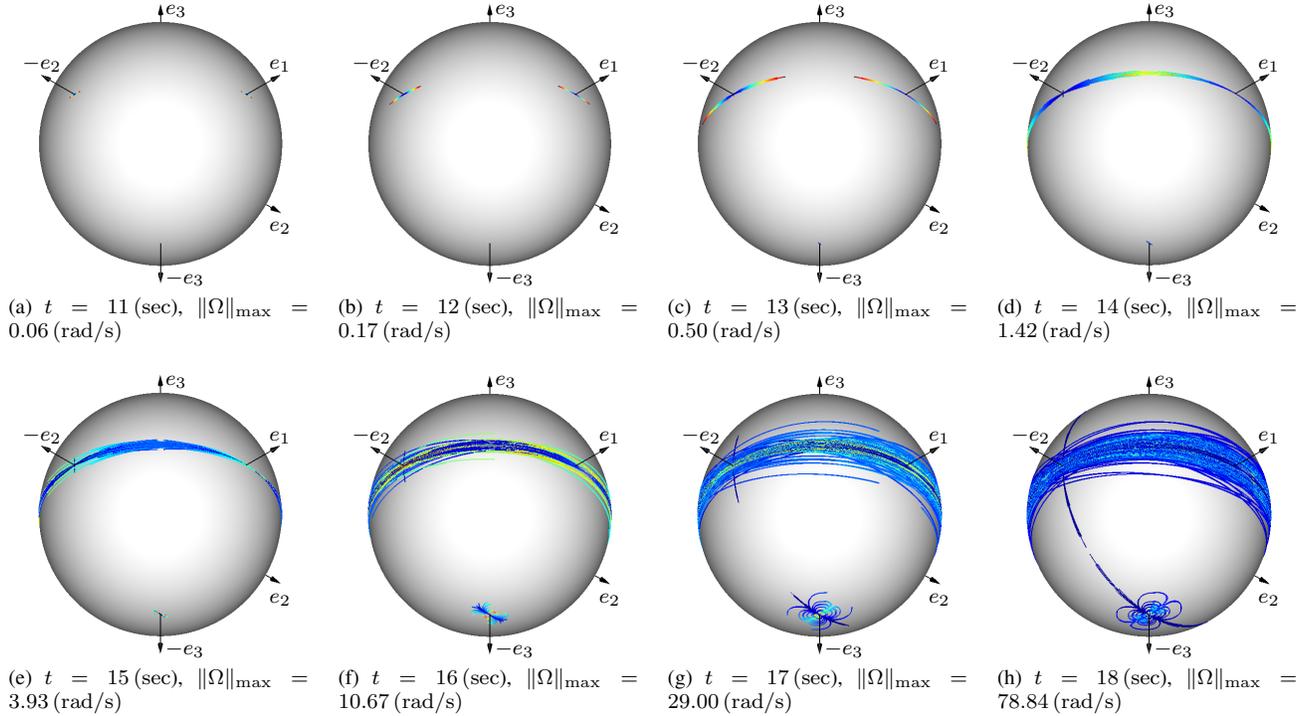

\subsubsection{Visualization of $W_s(\exp(\pi\hat e_1),0)$}

In~\cite{LeeLeoPICDC08}, a method to visualize a function or a trajectory on $\SO$ is proposed. Each column of a rotation matrix represents the direction of a body-fixed axis, and it evolves on $\Sph^2$. Therefore, a trajectory on $\SO$ can be visualized by three curves on a sphere, representing the trajectory of three columns of a rotation matrix. The direction of the angular velocity should be chosen such that the corresponding time-derivative of the rotation matrix is tangent to the curve, and the magnitude of angular velocity can be illustrated by color shading. An example of visualizing a rotation about a single axis is illustrated in \reffig{visSO3_demo}.

\setlength{\unitlength}{0.1\columnwidth}
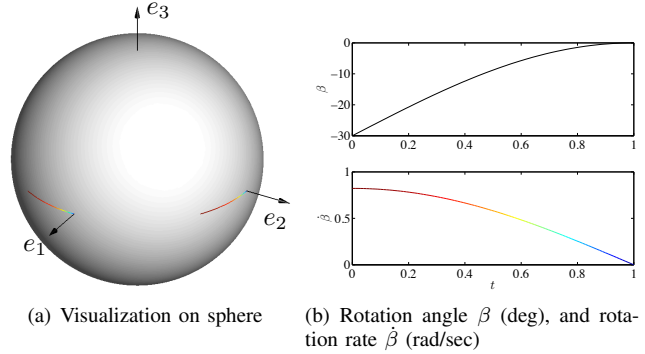
\begin{figure}
\centerline{
\subfigure[Visualization on sphere]{
\begin{picture}(4.5,4.5)(0,0)
\put(0,0){\includegraphics[width=0.45\columnwidth]{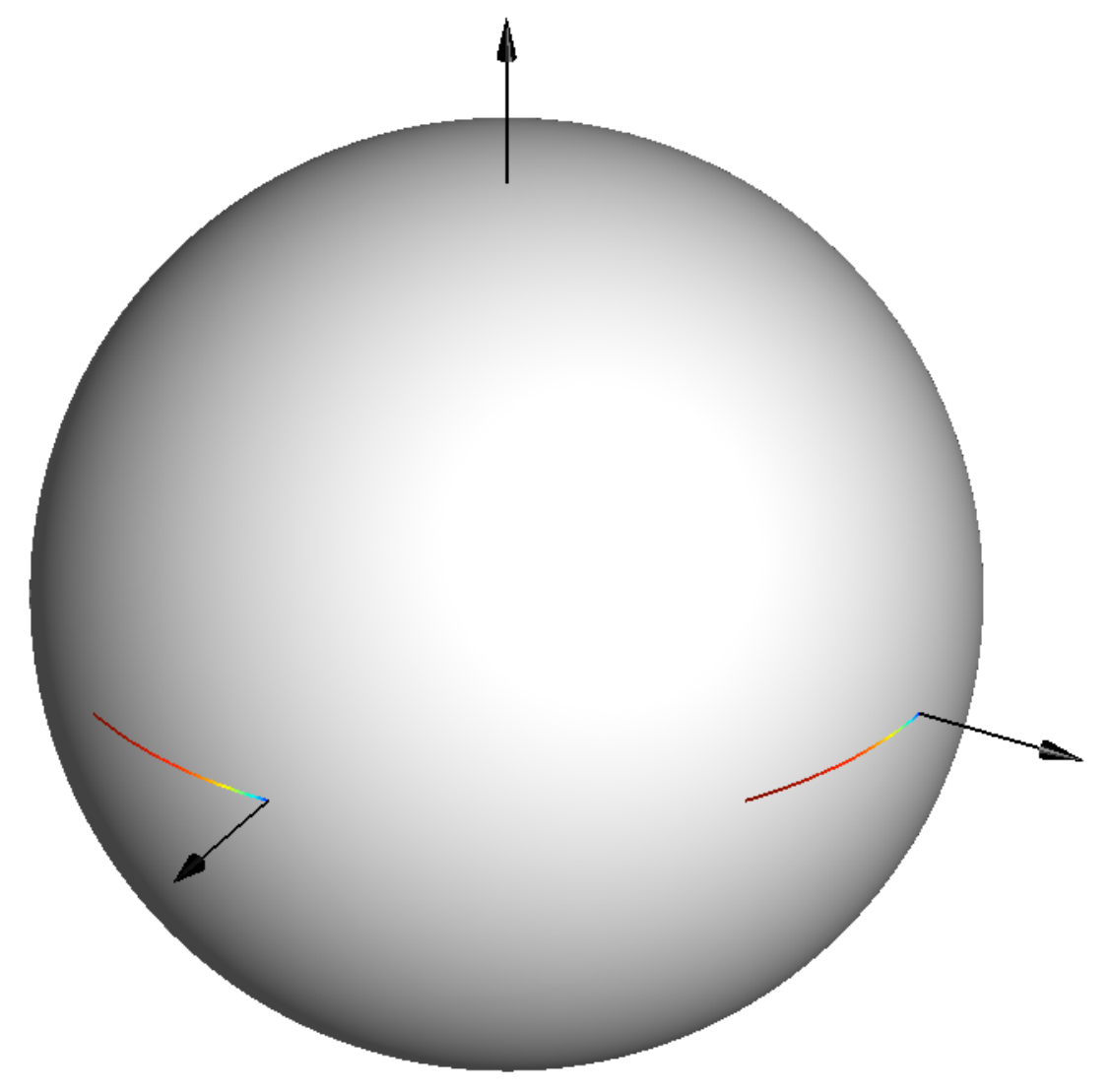}}
\put(4.0,1.0){$e_2$}\put(0.3,0.6){$e_1$}\put(2.2,4.25){$e_3$}
\end{picture}}
\subfigure[Rotation angle $\beta$ (deg), and rotation rate $\dot\beta$ (rad/sec)]{
\includegraphics[width=0.50\columnwidth]{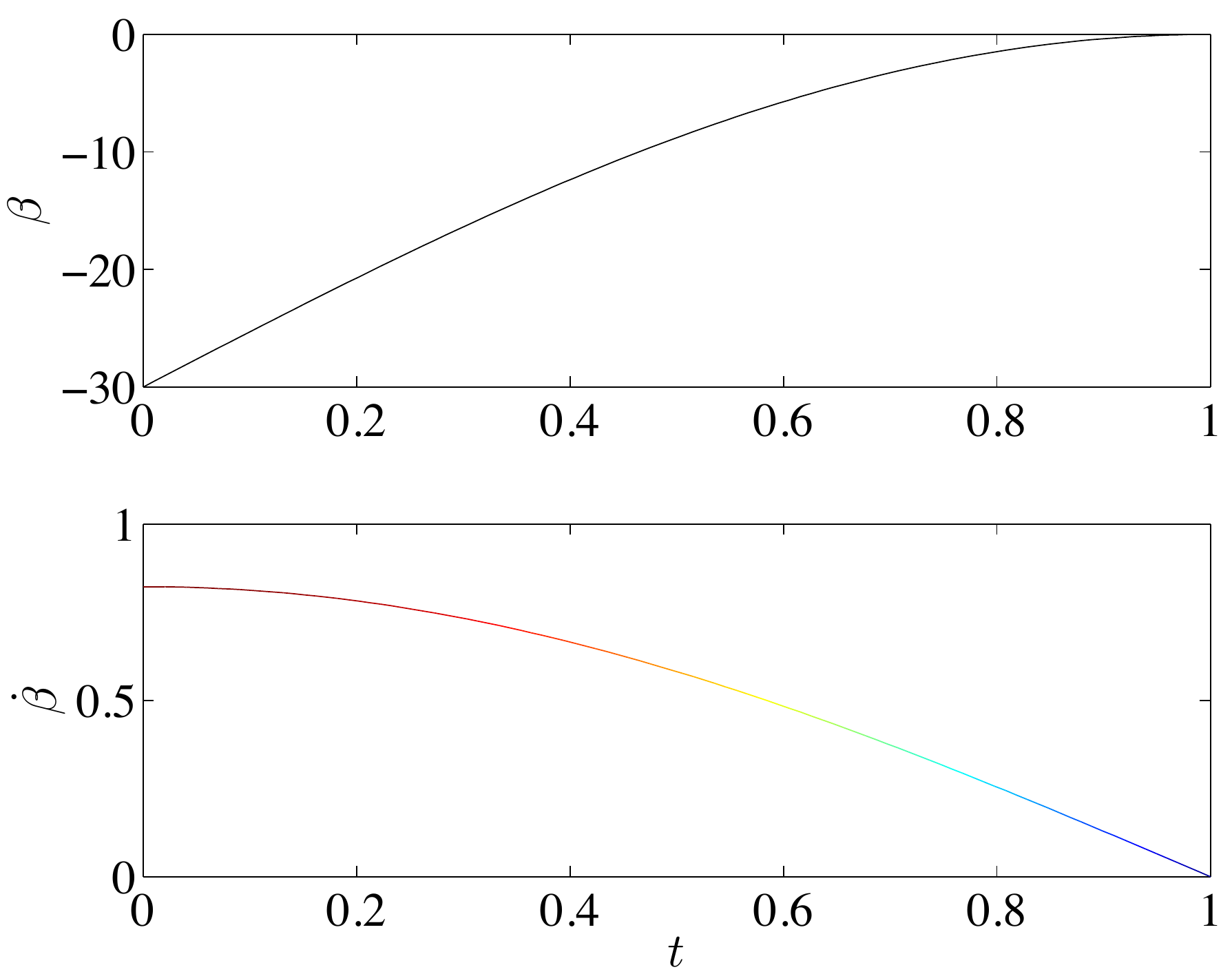}}
}
\caption{Visualization of an attitude maneuver: $R(t)=\exp(\beta(t)\hat e_3)$ for $0\leq t\leq 1$, where $\beta(t)=\frac{\pi}{6}(\sin\frac{\pi}{2}t -1)$. This maneuver corresponds to a rotation about the $e_3$ axis by $30^\circ$ to $R(1)=I$. The trajectory of the $i$-th column of $R(t)$ representing the direction of the $i$-th body-fixed axis is illustrated on a sphere for $i\in\{1,2,3\}$ (left). As the third body-fixed axis does not move during this maneuver, it is represented by a single point along the $e_3$ axis on the sphere. The direction of $\dot R(t)$ is tangent to these curves, and the magnitude of $\dot R(t)$ is denoted by color shading, according to the magnitude of the rotation rate (right).}\label{fig:visSO3_demo}
\end{figure}

We choose 112 points on the surface of $B_\delta\subset E^s_{loc} (\exp(\pi\hat e_1),0)$ with $\delta=10^{-6}$, and each point is integrated backward using \refeqn{qk}, \refeqn{wk} with timestep $h=0.002$. The resulting trajectories are illustrated in \reffig{SM1} for several values of $t$. 

In each figure, three body-fixed axes of the desired attitude $R_d=[e_1,e_2,e_3]$, and three body-fixed axes of the additional equilibrium attitude $\exp(\pi\hat e_1)=[e_1,-e_2,-e_3]$ are shown. From these computational results, we observe the following characteristics on the stable manifold $W_s(\exp(\pi\hat e_1),0)$:

\begin{itemize}
\item When $t\leq 15$, the trajectories in $W_s(\exp(\pi\hat e_1),0)$ are close to rotations about the third body-fixed axis $e_3$ to $\exp(\pi\hat e_1)$. This is consistent with the linearized dynamics, where the eigenvalue of the third mode, corresponding to the rotations about $e_3$, has the fastest convergence rate, as seen in \refeqn{v3_SO31}. 
\item When $t\geq 15$, the first mode representing the rotations about $e_1$ starts to appear, followed by the second mode representing the rotation about $e_2$. This corresponds to the fact that the first mode has a faster convergence rate than the second mode, i.e. $|\lambda_1|>|\lambda_2|$. 
\item As $t$ is increased further, the third body-fixed axis leaves the neighborhood of $-e_3$, and it exhibit the following pattern:

\centerline{\includegraphics[width=0.32\columnwidth]{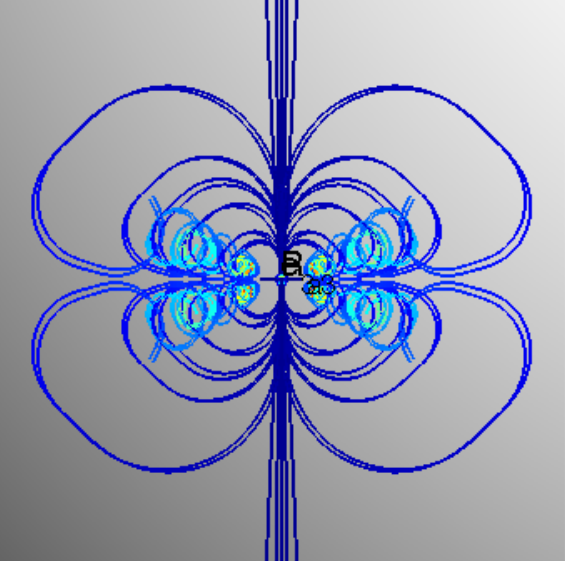}}
\item The stable manifold $W_s(\exp(\pi\hat e_1),0)$ covers a certain part of $\SO$, when projected on to it. So, when an initial attitude is chosen such that its third body-fixed axis is sufficiently close to $-e_3$, there  possibly exist multiple initial angular velocities such that the corresponding solution converges to $\exp(\pi\hat e_1)$ instead of the desired attitude $R_d=I$.
\end{itemize}

\subsubsection{Visualization of $W_s(\exp(\pi\hat e_2),0)$}

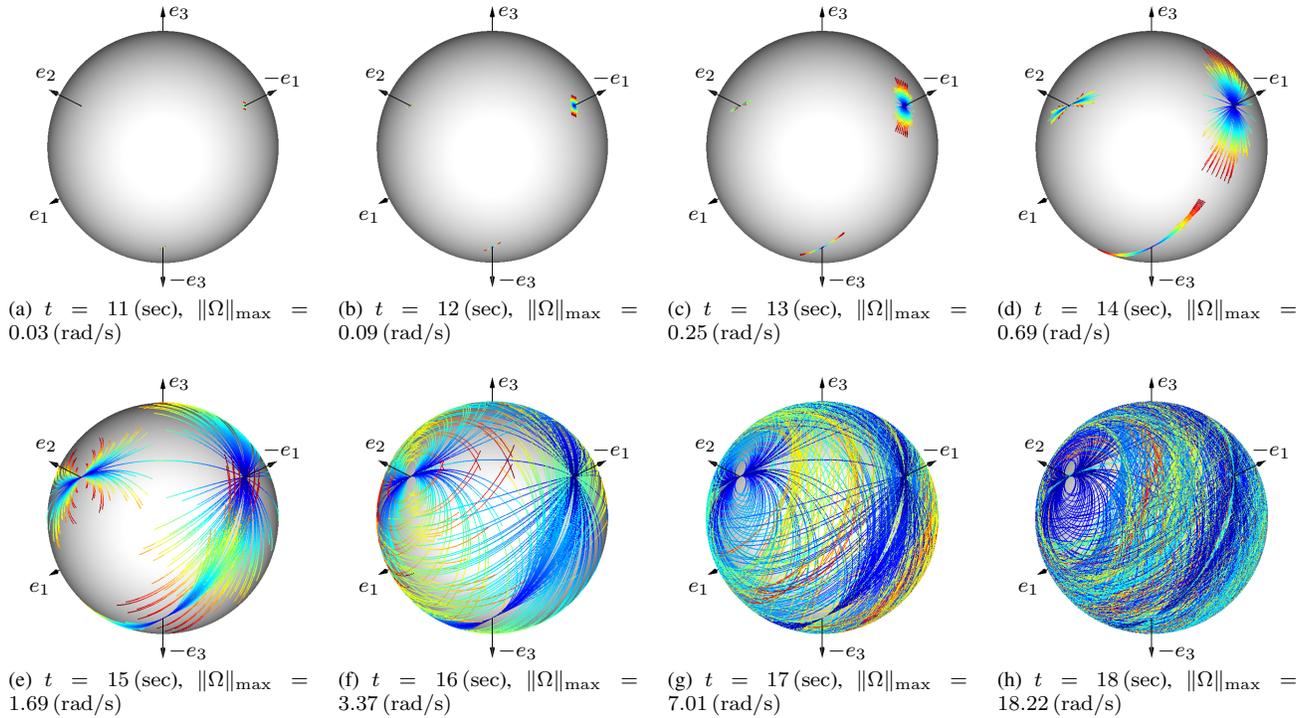
\begin{figure*}
\footnotesize\selectfont
\centerline{
\subfigure[$t=11\,(\mathrm{sec}$), $\|\Omega\|_{\max}=0.03\,(\mathrm{rad/s})$]{
\begin{picture}(4.5,4.5)(0,0)
\put(0,0){\includegraphics[width=0.45\columnwidth]{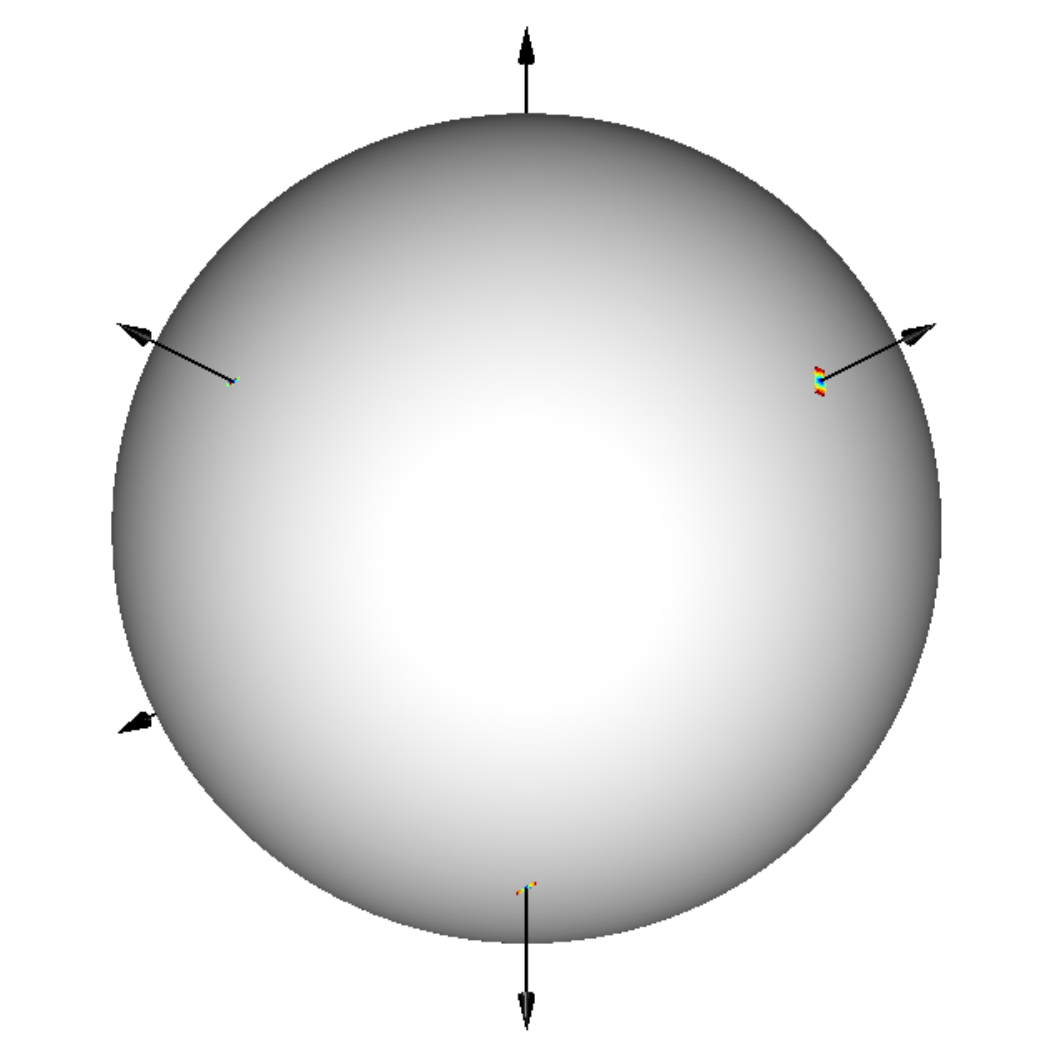}}
\put(0.2,1.1){$e_1$}\put(0.3,3.3){$e_2$}\put(2.35,4.25){$e_3$}
\put(3.8,3.2){$-e_1$}\put(2.35,0.1){$-e_3$}
\end{picture}}
\hspace*{0.2cm}
\subfigure[$t=12\,(\mathrm{sec}$), $\|\Omega\|_{\max}=0.09\,(\mathrm{rad/s})$]{
\begin{picture}(4.5,4.5)(0,0)
\put(0,0){\includegraphics[width=0.45\columnwidth]{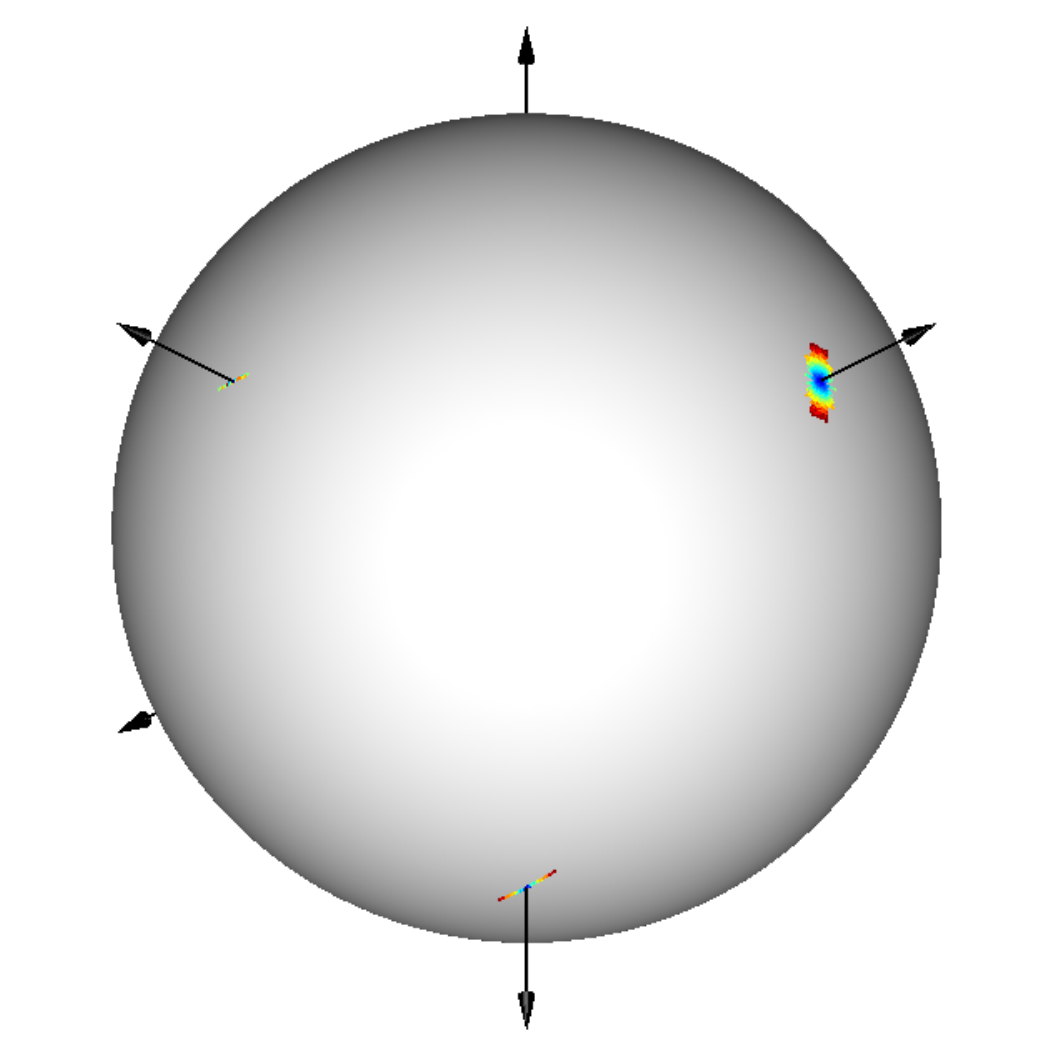}}
\put(0.2,1.1){$e_1$}\put(0.3,3.3){$e_2$}\put(2.35,4.25){$e_3$}
\put(3.8,3.2){$-e_1$}\put(2.35,0.1){$-e_3$}
\end{picture}}
\hspace*{0.2cm}
\subfigure[$t=13\,(\mathrm{sec}$), $\|\Omega\|_{\max}=0.25\,(\mathrm{rad/s})$]{
\begin{picture}(4.5,4.5)(0,0)
\put(0,0){\includegraphics[width=0.45\columnwidth]{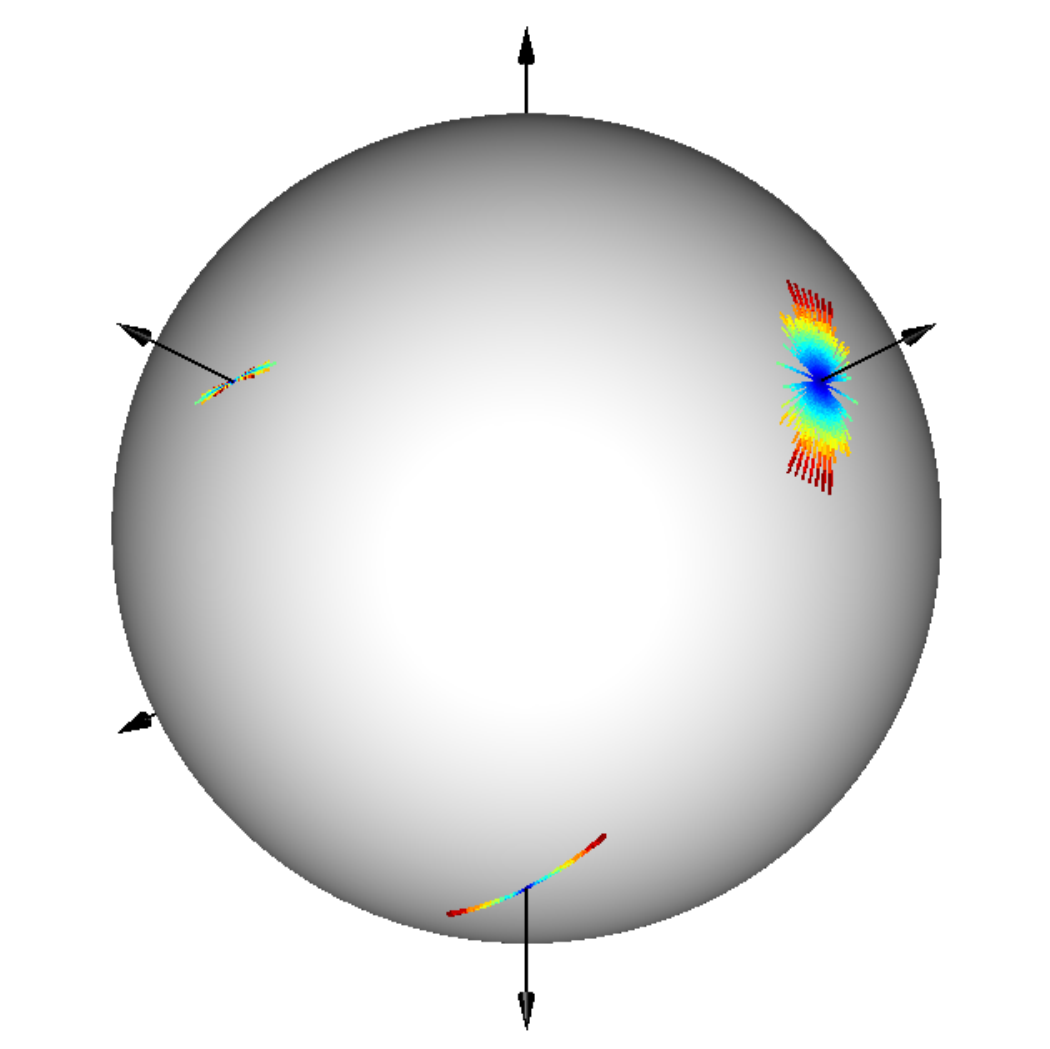}}
\put(0.2,1.1){$e_1$}\put(0.3,3.3){$e_2$}\put(2.35,4.25){$e_3$}
\put(3.8,3.2){$-e_1$}\put(2.35,0.1){$-e_3$}
\end{picture}}
\hspace*{0.2cm}
\subfigure[$t=14\,(\mathrm{sec}$), $\|\Omega\|_{\max}=0.69\,(\mathrm{rad/s})$]{
\begin{picture}(4.5,4.5)(0,0)
\put(0,0){\includegraphics[width=0.45\columnwidth]{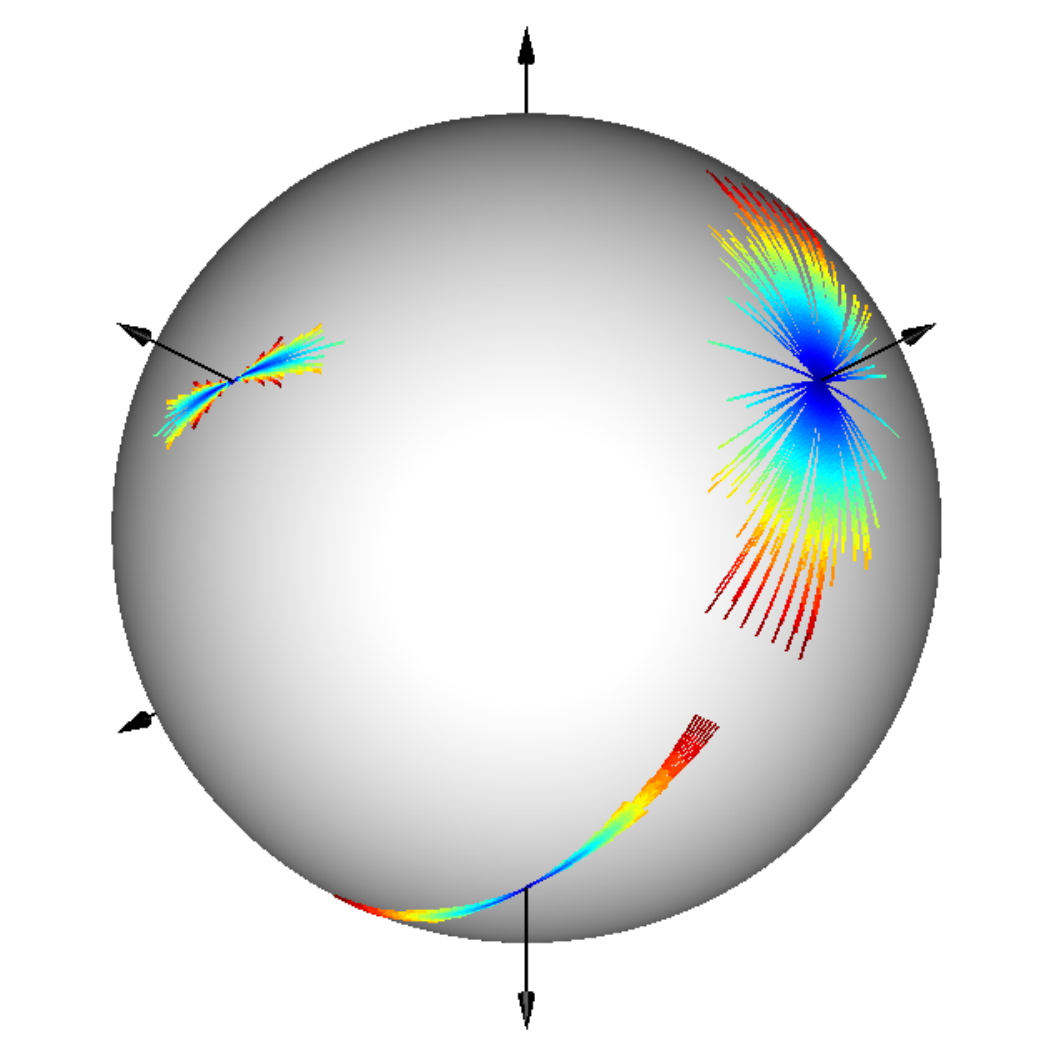}}
\put(0.2,1.1){$e_1$}\put(0.3,3.3){$e_2$}\put(2.35,4.25){$e_3$}
\put(3.8,3.2){$-e_1$}\put(2.35,0.1){$-e_3$}
\end{picture}}
}
\centerline{
\subfigure[$t=15\,(\mathrm{sec}$), $\|\Omega\|_{\max}=1.69\,(\mathrm{rad/s})$]{
\begin{picture}(4.5,4.5)(0,0)
\put(0,0){\includegraphics[width=0.45\columnwidth]{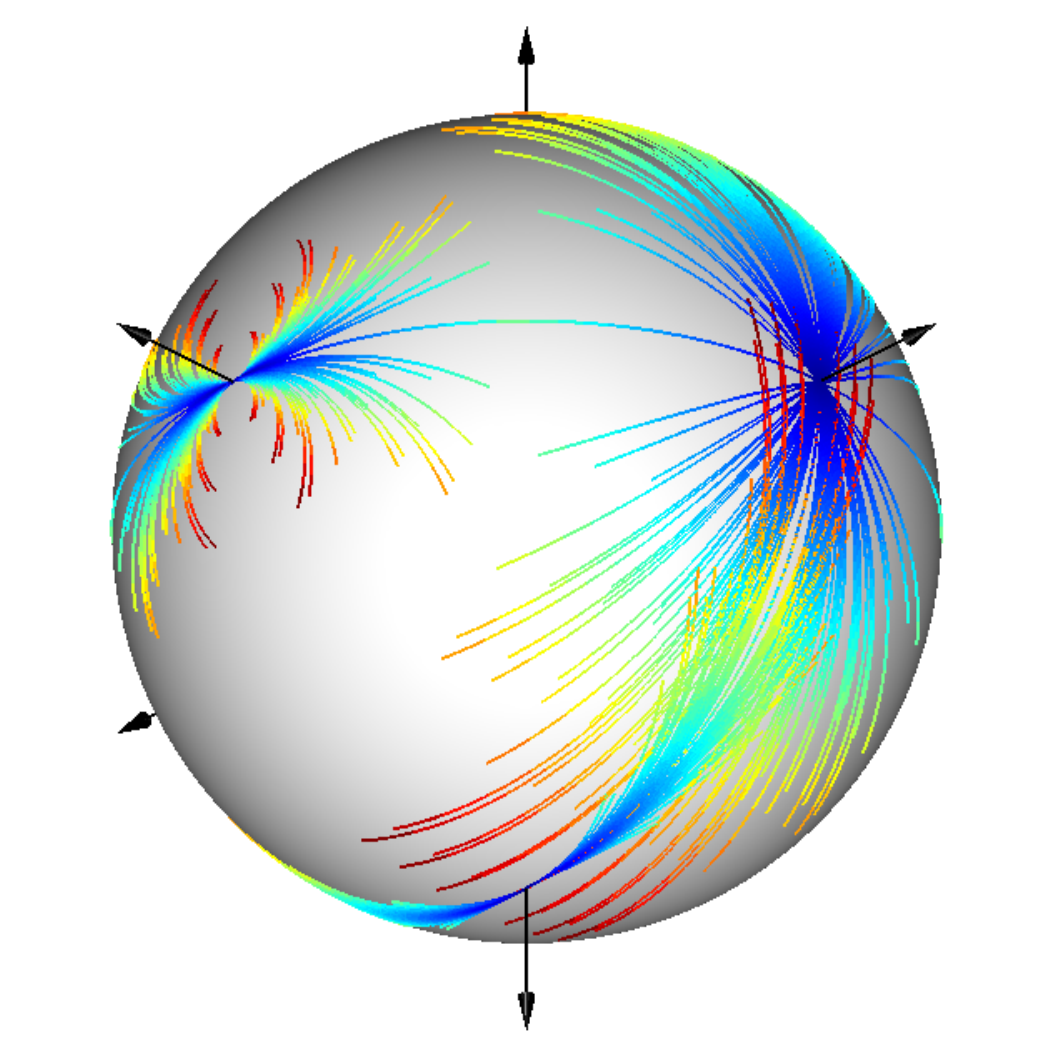}}
\put(0.2,1.1){$e_1$}\put(0.3,3.3){$e_2$}\put(2.35,4.25){$e_3$}
\put(3.8,3.2){$-e_1$}\put(2.35,0.1){$-e_3$}
\end{picture}}
\hspace*{0.2cm}
\subfigure[$t=16\,(\mathrm{sec}$), $\|\Omega\|_{\max}=3.37\,(\mathrm{rad/s})$]{
\begin{picture}(4.5,4.5)(0,0)
\put(0,0){\includegraphics[width=0.45\columnwidth]{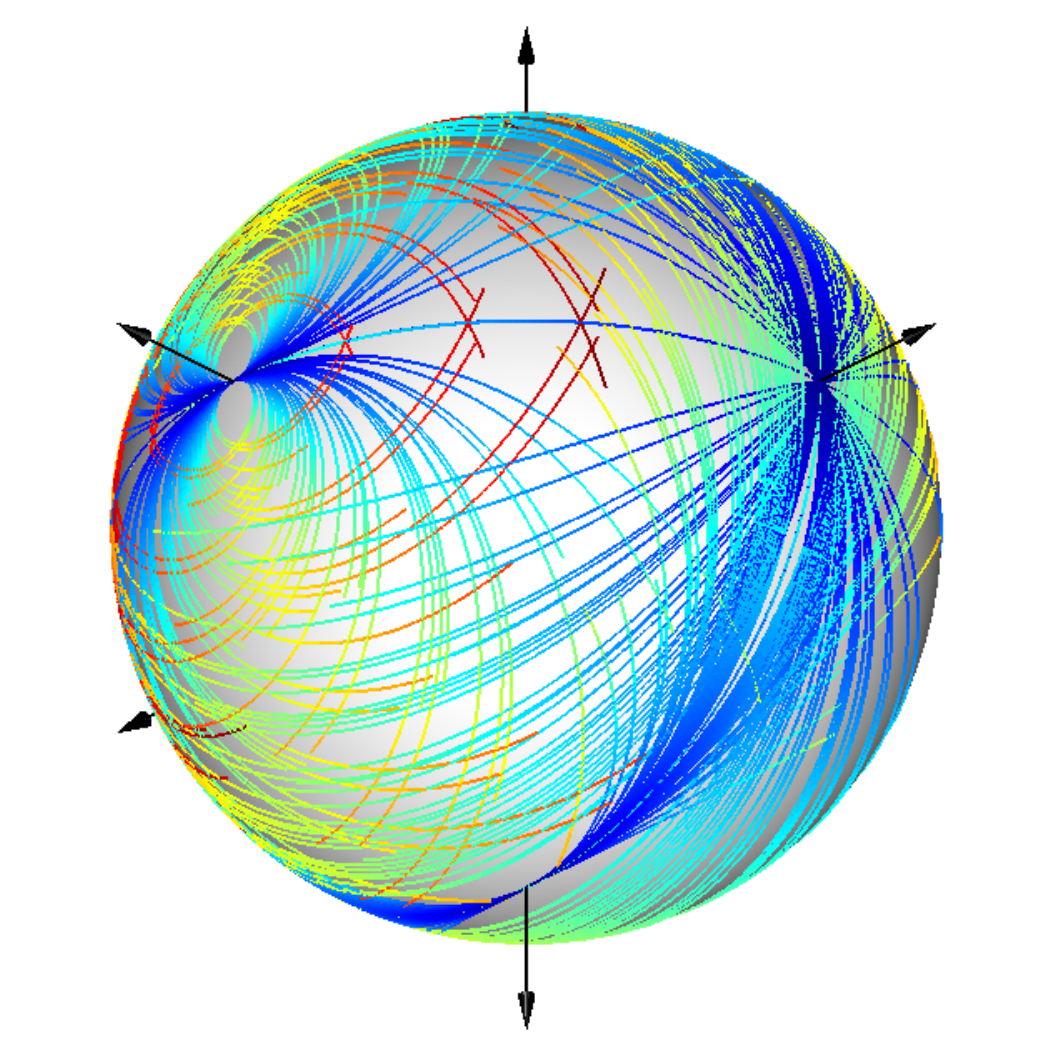}}
\put(0.2,1.1){$e_1$}\put(0.3,3.3){$e_2$}\put(2.35,4.25){$e_3$}
\put(3.8,3.2){$-e_1$}\put(2.35,0.1){$-e_3$}
\end{picture}}
\hspace*{0.2cm}
\subfigure[$t=17\,(\mathrm{sec}$), $\|\Omega\|_{\max}=7.01\,(\mathrm{rad/s})$]{
\begin{picture}(4.5,4.5)(0,0)
\put(0,0){\includegraphics[width=0.45\columnwidth]{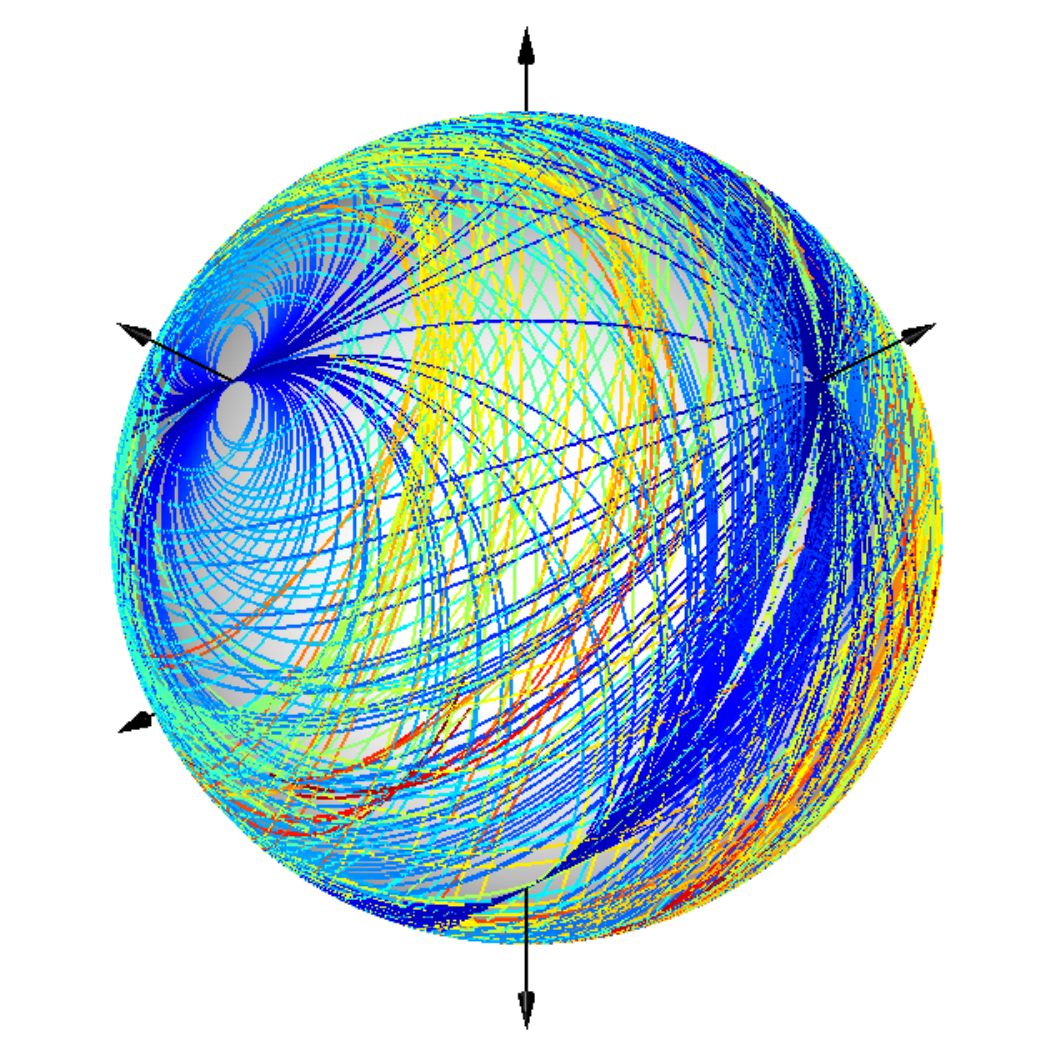}}
\put(0.2,1.1){$e_1$}\put(0.3,3.3){$e_2$}\put(2.35,4.25){$e_3$}
\put(3.8,3.2){$-e_1$}\put(2.35,0.1){$-e_3$}
\end{picture}}
\hspace*{0.2cm}
\subfigure[$t=18\,(\mathrm{sec}$), $\|\Omega\|_{\max}=18.22\,(\mathrm{rad/s})$]{
\begin{picture}(4.5,4.5)(0,0)
\put(0,0){\includegraphics[width=0.45\columnwidth]{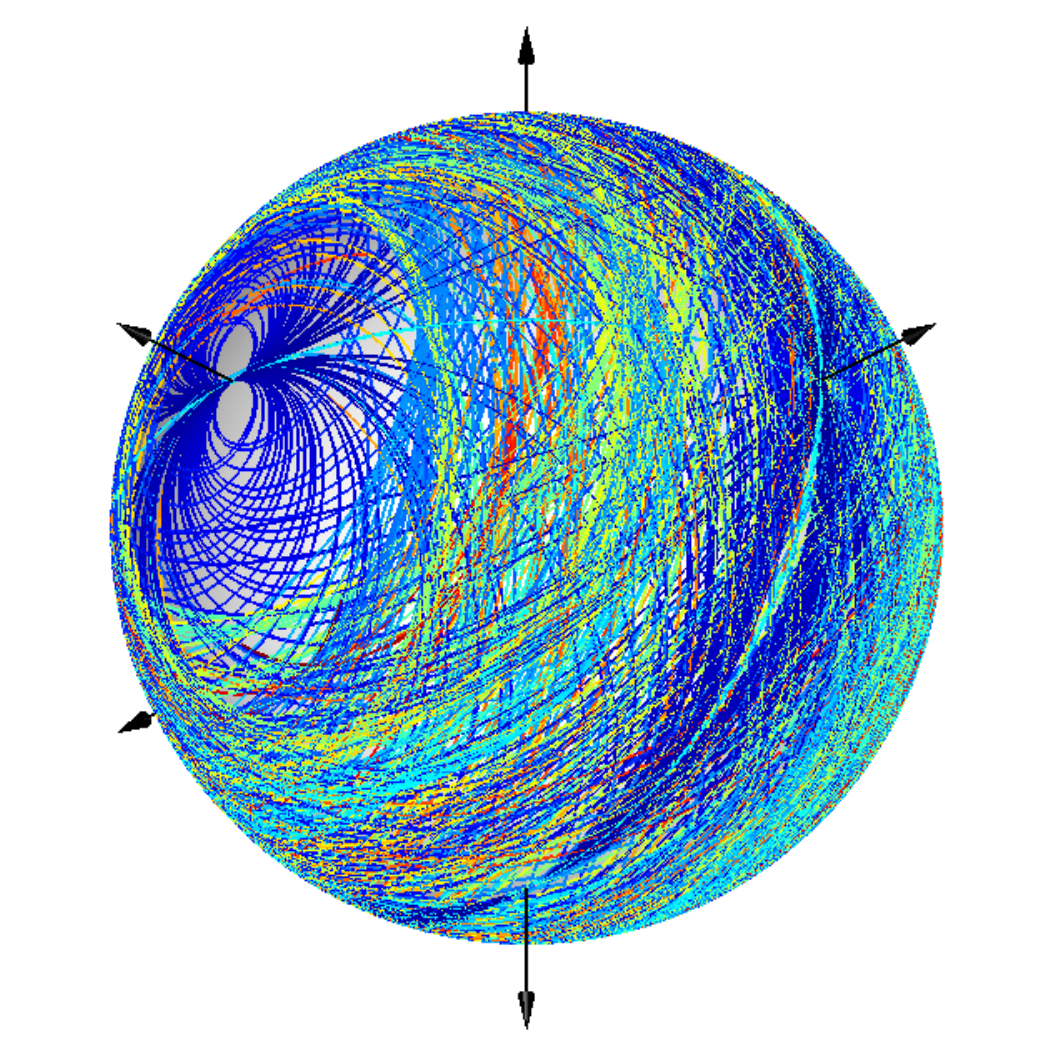}}
\put(0.2,1.1){$e_1$}\put(0.3,3.3){$e_2$}\put(2.35,4.25){$e_3$}
\put(3.8,3.2){$-e_1$}\put(2.35,0.1){$-e_3$}
\end{picture}}
}
\caption{Stable manifold to $(\exp(\pi\hat e_2),0)=([-e_1,e_2,-e_3],0)$ represented by $\{ \mathcal{F}^{-t} (B_\delta)\}_{t>0}$ with $\delta=10^{-6}$ for several values of $t$.}\label{fig:SM2}
\end{figure*}

We choose 544 points on the surface of $B_\delta\subset E^s_{loc} (\exp(\pi\hat e_2),0)$ with $\delta=10^{-6}$, and each point is integrated backward using \refeqn{qk}, \refeqn{wk} with timestep $h=0.002$. The resulting trajectories are illustrated in \reffig{SM2} for several values of $t$. 

In each figure, three body-fixed axes of the desired attitude $R_d=[e_1,e_2,e_3]$, and three body-fixed axes of the additional equilibrium attitude $\exp(\pi\hat e_2)=[-e_1,e_2,-e_3]$ are shown. From these computational results, we observe the following characteristics on the stable manifold $W_s(\exp(\pi\hat e_2),0)$:

\begin{itemize}
\item When $t\leq 12$, the trajectories in $W_s(\exp(\pi\hat e_2),0)$ is close to the rotations about the second body-fixed axis $e_2$. As $t$ increases, rotations about $e_3$ starts to appear. This corresponds to the linearized dynamics where the second mode representing rotations about $e_2$ has the fastest convergence rate, followed by the third mode at \refeqn{v2_SO32}.
\item As $t$ is increased further, nonlinear modes become dominant. The trajectories in $W_s(\exp(\pi\hat e_2),0)$ almost cover \SO. This suggests that for any initial attitude, we can choose several initial angular velocities such that the corresponding solutions converges to $\exp(\pi\hat e_2)$.
\end{itemize}

\subsubsection{Visualization of $W_s(\exp(\pi\hat e_3),0)$}

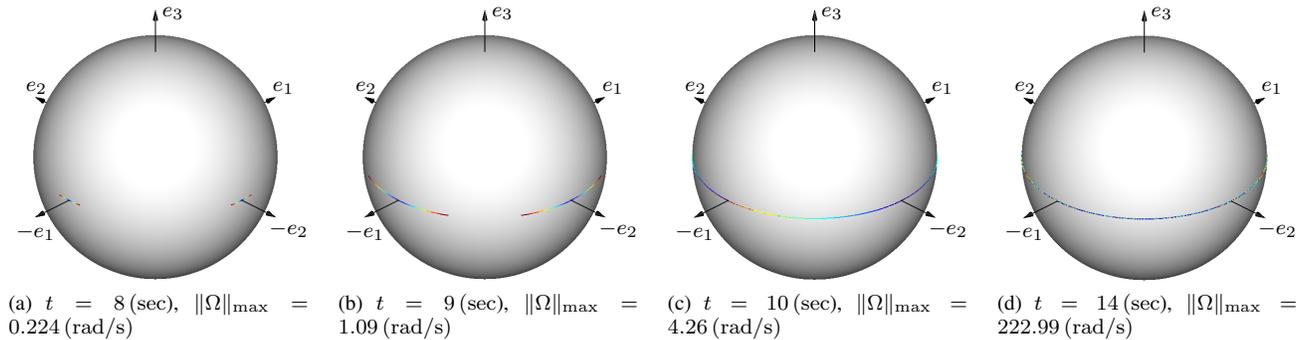
\begin{figure*}
\footnotesize\selectfont
\centerline{
\subfigure[$t=8\,(\mathrm{sec}$), $\|\Omega\|_{\max}=0.224\,(\mathrm{rad/s})$]{
\begin{picture}(4.5,4.5)(0,0)
\put(0,0){\includegraphics[width=0.43\columnwidth]{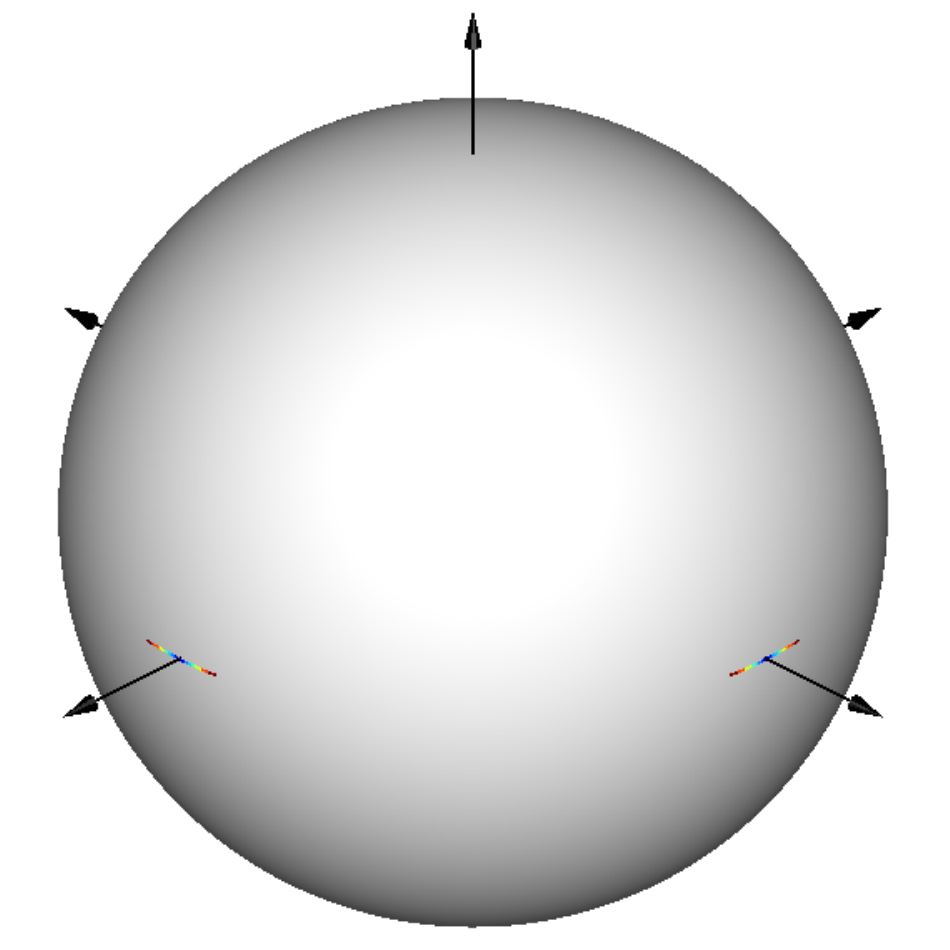}}
\put(0.0,0.75){$-e_1$}\put(0.15,3.0){$e_2$}\put(2.25,4.15){$e_3$}
\put(3.95,3.0){$e_1$}\put(3.90,0.8){$-e_2$}
\end{picture}}
\hspace*{0.2cm}
\subfigure[$t=9\,(\mathrm{sec}$), $\|\Omega\|_{\max}=1.09\,(\mathrm{rad/s})$]{
\begin{picture}(4.5,4.5)(0,0)
\put(0,0){\includegraphics[width=0.43\columnwidth]{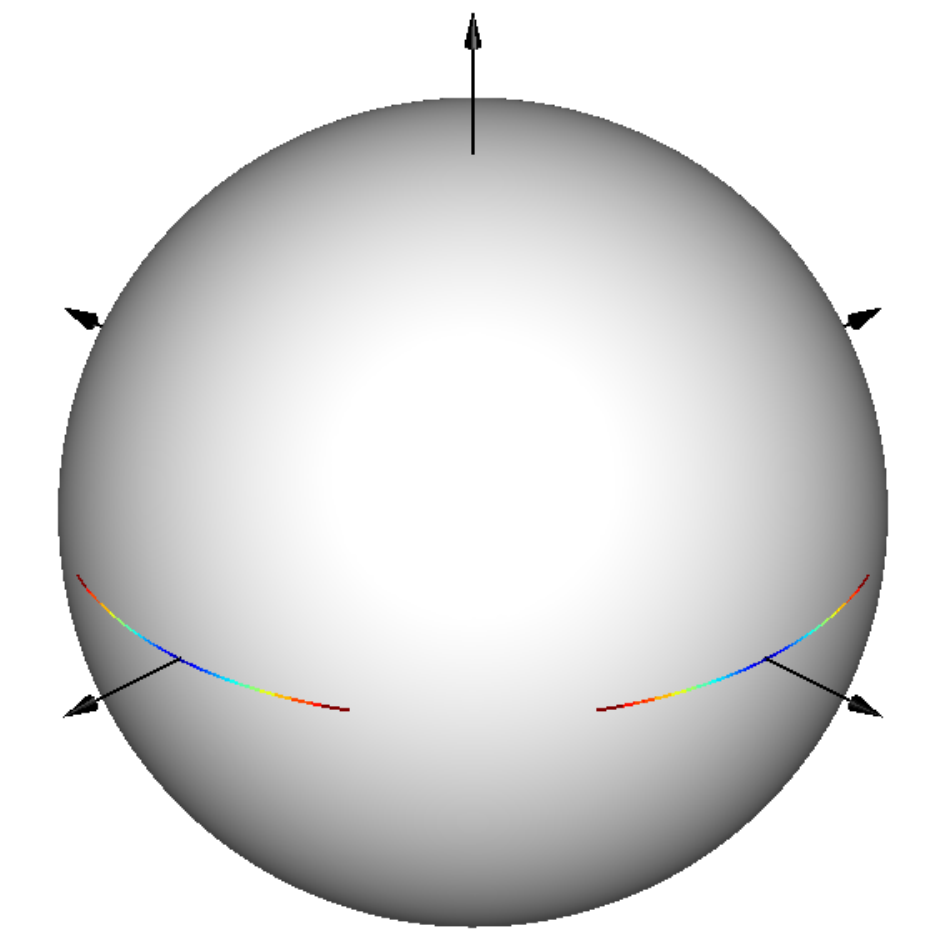}}
\put(0.0,0.75){$-e_1$}\put(0.15,3.0){$e_2$}\put(2.25,4.15){$e_3$}
\put(3.95,3.0){$e_1$}\put(3.90,0.8){$-e_2$}
\end{picture}}
\hspace*{0.2cm}
\subfigure[$t=10\,(\mathrm{sec}$), $\|\Omega\|_{\max}=4.26\,(\mathrm{rad/s})$]{
\begin{picture}(4.5,4.5)(0,0)
\put(0,0){\includegraphics[width=0.43\columnwidth]{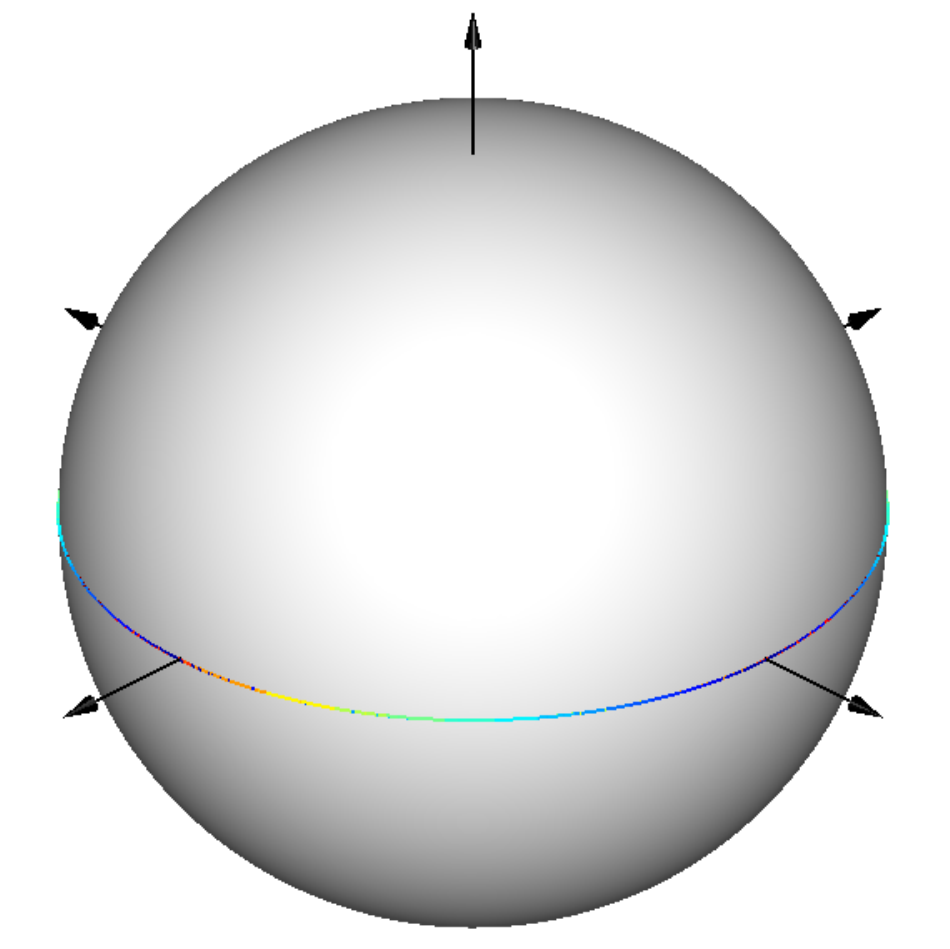}}
\put(0.0,0.75){$-e_1$}\put(0.15,3.0){$e_2$}\put(2.25,4.15){$e_3$}
\put(3.95,3.0){$e_1$}\put(3.90,0.8){$-e_2$}
\end{picture}}
\hspace*{0.2cm}
\subfigure[$t=14\,(\mathrm{sec}$), $\|\Omega\|_{\max}=222.99\,(\mathrm{rad/s})$]{
\begin{picture}(4.5,4.5)(0,0)
\put(0,0){\includegraphics[width=0.43\columnwidth]{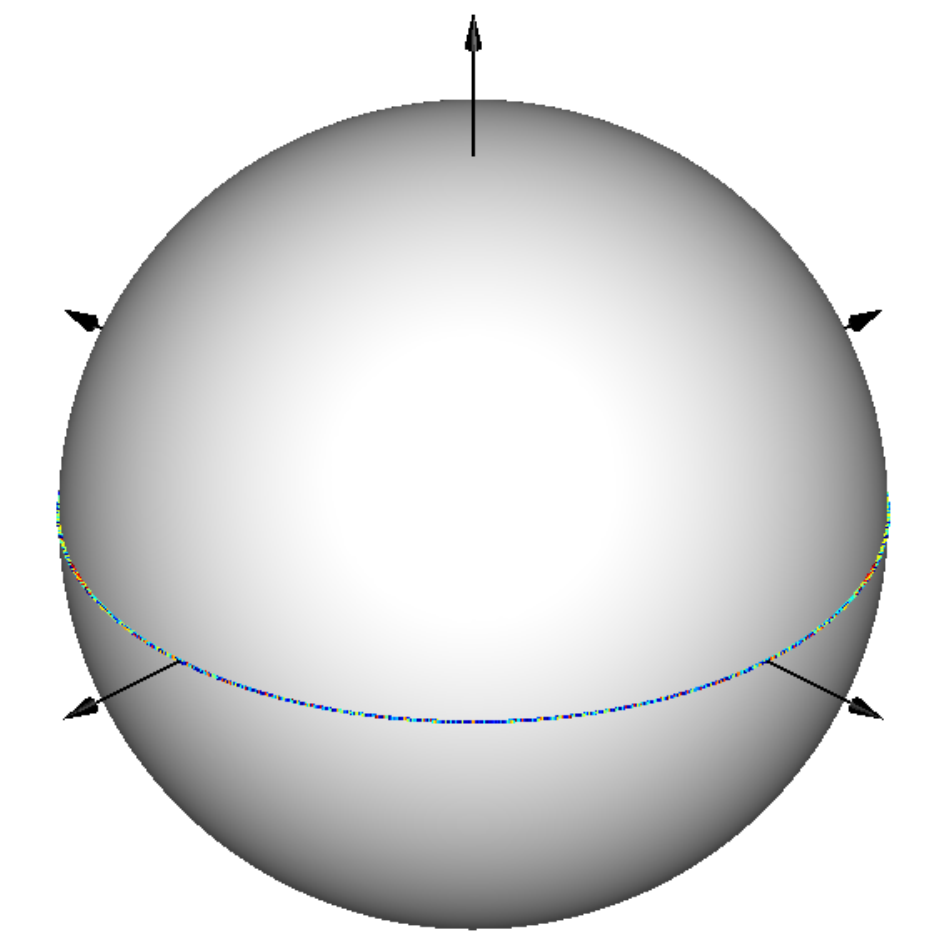}}
\put(0.0,0.75){$-e_1$}\put(0.15,3.0){$e_2$}\put(2.25,4.15){$e_3$}
\put(3.95,3.0){$e_1$}\put(3.90,0.8){$-e_2$}
\end{picture}}
}

\caption{Stable manifold to $(\exp(\pi\hat e_3),0)=([-e_1,-e_2,e_3],0)$ represented by $\{ \mathcal{F}^{-t} (B_\delta)\}_{t>0}$ with $\delta=10^{-6}$ for several values of $t$.}\label{fig:SM3}
\end{figure*}

Similarly, we choose 976 points on the surface of $B_\delta\subset E^s_{loc} (\exp(\pi\hat e_3),0)$ with $\delta=10^{-6}$, and each point is integrated backward using \refeqn{qk}, \refeqn{wk} with timestep $h=0.002$. The resulting trajectories are illustrated in \reffig{SM3} for several values of $t$. 

At each figure, three body-fixed axes of the desired attitude $R_d=[e_1,e_2,e_3]$, and three body-fixed axes of the additional equilibrium attitude $\exp(\pi\hat e_3)=[-e_1,-e_2,e_3]$ are shown. From these computational results, we observe the following characteristics on the stable manifold $W_s(\exp(\pi\hat e_3),0)$:

\begin{itemize}
\item When $t\leq 8$, the trajectories in $W_s(\exp(\pi\hat e_3),0)$ are close to the rotations about the third body-fixed axis $e_3$. This corresponds to the linearized dynamics where the fifth mode representing rotations about $e_3$ has the fastest convergence rate given in \refeqn{v5_SO33}.
\item The rotations about $e_3$ are still dominant, even as $t$ is increased further. For the given simulation times, all trajectories in $W_s(\exp(\pi\hat e_3),0)$ are close to rotations about $e_3$.
\end{itemize}

\section{Conclusions}

Stable manifolds of saddle points that arise in the closed-loop dynamics of two pendulum models are characterized numerically, and several properties are observed. Although the analytical and computational results have been presented for a spherical pendulum and a 3D pendulum, the methods presented naturally extend to any closed loop attitude control system with configurations in either $\Sph^2$ or $\SO$.    

\bibliography{CDC11}

\begin{thebibliography}{10}
\providecommand{\url}[1]{#1}
\csname url@rmstyle\endcsname
\providecommand{\newblock}{\relax}
\providecommand{\bibinfo}[2]{#2}
\providecommand\BIBentrySTDinterwordspacing{\spaceskip=0pt\relax}
\providecommand\BIBentryALTinterwordstretchfactor{4}
\providecommand\BIBentryALTinterwordspacing{\spaceskip=\fontdimen2\font plus
\BIBentryALTinterwordstretchfactor\fontdimen3\font minus
  \fontdimen4\font\relax}
\providecommand\BIBforeignlanguage[2]{{%
\expandafter\ifx\csname l@#1\endcsname\relax
\typeout{** WARNING: IEEEtran.bst: No hyphenation pattern has been}%
\typeout{** loaded for the language `#1'. Using the pattern for}%
\typeout{** the default language instead.}%
\else
\language=\csname l@#1\endcsname
\fi
#2}}

\bibitem{ChSaMcCSM11}
N.~A. Chaturvedi, A.~K. Sanyal, and N.~H. McClamroch, ``Rigid body attitude
  control: Using rotation matrices for continuous, singularity-free control
  laws,'' \emph{IEEE Control Systems Magazine}, p. accepted, 2011.

\bibitem{ChMcIJRNC07}
N.~A. Chaturvedi and N.~H. McClamroch, ``Asymptotic stabilization of the
  hanging equilibrium manifold of the 3{D} pendulum,'' \emph{International
  Journal of Robust and Nonlinear Control}, pp. 1435--1454, 2007.

\bibitem{ChMcBeAut08}
N.~A. Chaturvedi, N.~H. McClamroch, and D.~S. Bernstein, ``Stabilization of a
  3{D} axially symmetric pendulum,'' \emph{Automatica}, pp. 2258--2265, 2008.

\bibitem{ChaMcCITAC09}
N.~Chaturvedi, N.~H. McClamroch, and D.~Bernstein, ``Asymptotic smooth
  stabilization of the inverted 3-{D} pendulum,'' \emph{IEEE Transactions on
  Automatic Control}, vol.~54, no.~6, pp. 1204--1215, 2009.

\bibitem{BulMurN95}
F.~Bullo, R.~M. Murray, and A.~Sarti, ``Control on the sphere and reduced
  attitude stabilization,'' in \emph{IFAC Symposium on Nonlinear Control
  Systems}, vol.~2, 1995, pp. 495--501.

\bibitem{BulLew05}
F.~Bullo and A.~Lewis, \emph{Geometric control of mechanical systems}.\hskip
  1em plus 0.5em minus 0.4em\relax Springer, 2005.

\bibitem{LeeLeoIJNME08}
T.~Lee, M.~Leok, and N.~H. McClamroch, ``Lagrangian mechanics and variational
  integrators on two-spheres,'' \emph{International Journal for Numerical
  Methods in Engineering}, vol.~79, no.~9, pp. 1147--1174, 2009.

\bibitem{Kuz98}
Y.~Kuznetsov, \emph{Elements of Applied Bifurcation Theory}.\hskip 1em plus
  0.5em minus 0.4em\relax Springer, 1998.

\bibitem{KraOsiIJBC05}
B.~Krauskopf, H.~Osinga, E.~Doedel, M.~Henderson, J.~Guckenheimer,
  A.~Vladimirsky, M.~Dellnitz, and O.~Junge, ``A survey of methods for
  computing (un)stable manifolds of vector fields a survey of methods for
  computing (un)stable manifolds of vector fields,'' \emph{International
  Journal of Bifurcation and Chaos}, vol.~15, no.~3, pp. 763--791, 2005.

\bibitem{HaiLub00}
E.~Hairer, C.~Lubich, and G.~Wanner, \emph{Geometric numerical integration},
  ser. Springer Series in Computational Mechanics 31.\hskip 1em plus 0.5em
  minus 0.4em\relax Springer, 2000.

\bibitem{MarWesAN01}
J.~Marsden and M.~West, ``Discrete mechanics and variational integrators,'' in
  \emph{Acta Numerica}.\hskip 1em plus 0.5em minus 0.4em\relax Cambridge, 2001,
  vol.~10, pp. 317--514.

\bibitem{ChaLeeJNS11}
N.~Chaturvedi, T.~Lee, M.~Leok, and N.~H. McClamroch, ``Nonlinear dynamics of
  the 3{D} pendulum,'' \emph{Journal of Nonlinear Science}, vol.~21, no.~1, pp.
  3--21, 2011.

\bibitem{LeePACC11}
T.~Lee, ``Geometric tracking control of the attitude dynamics of a rigid body
  on {SO(3)},'' in \emph{Proceeding of the American Control Conference}, 2011,
  accepted.

\bibitem{BhaBerSCL00}
S.~Bhat and D.~Bernstein, ``A topological obstruction to continuous global
  stabilization of rotational motion and the unwinding phenomenon,''
  \emph{Systems and Control Letters}, vol.~39, no.~1, pp. 66--73, 2000.

\bibitem{KodPICDC98}
D.~Koditschek, ``Application of a new lyapunov function to global adaptive
  tracking,'' in \emph{Proceedings of the IEEE Conference on Decision and
  Control}, 1998, pp. 63--68.

\bibitem{LeeLeoPICCA05}
T.~Lee, M.~Leok, and N.~H. McClamroch, ``A {L}ie group variational integrator
  for the attitude dynamics of a rigid body with application to the 3{D}
  pendulum,'' in \emph{Proceedings of the {IEEE} {C}onference on {C}ontrol
  {A}pplication}, 2005, pp. 962--967.

\bibitem{LeeLeoCMDA07}
------, ``Lie group variational integrators for the full body problem in
  orbital mechanics,'' \emph{Celestial Mechanics and Dynamical Astronomy},
  vol.~98, no.~2, pp. 121--144, June 2007.

\bibitem{LeeLeoPICDC08}
------, ``Global symplectic uncertainty propagation on {S}{O}(3),'' in
  \emph{Proceedings of the IEEE Conference on Decision and Control}, 2008, pp.
  61--66.

\end{thebibliography}
\bibliographystyle{IEEEtran}

\end{document}